\documentclass[conference]{IEEEtran}
\ifCLASSINFOpdf
\else
\fi

\IEEEoverridecommandlockouts

\def\BibTeX{{\rm B\kern-.05em{\sc i\kern-.025em b}\kern-.08em
    T\kern-.1667em\lower.7ex\hbox{E}\kern-.125emX}}

\usepackage{enumerate}
\usepackage{amsmath}
\usepackage{mathrsfs}
\usepackage[dvipsnames]{xcolor}	
\usepackage{bm}
\usepackage{graphicx}
\usepackage{subcaption}
\usepackage{adjustbox}
\usepackage{soul}
\usepackage{listings}
\usepackage{array}
\usepackage{makecell}
\usepackage{mathtools}
\usepackage{abraces}
\usepackage{booktabs}
\usepackage{algorithm} 
\usepackage{algpseudocode} 
\usepackage{semantic}
\usepackage{syntax}
\usepackage{tikz}
\usepackage{pgfplots}
\usepackage{subcaption}
\usepackage{xspace}
\usepackage[frozencache,cachedir=minted-cache]{minted}
\usepackage{multirow}
\usepackage{paralist}
\usepackage[inline]{enumitem}
\usepackage{stmaryrd}
\usepackage{colortbl}

\usepackage[english]{babel}

\usepackage{hyperref}
\usepackage{cleveref}

\ifCLASSOPTIONcompsoc
  \usepackage[nocompress]{cite}
\else
  \usepackage{cite}
\fi
\ifCLASSINFOpdf
\else
\fi

\AtBeginDocument{\DeclareCaptionSubType{lstlisting}}
\crefname{sublstlisting}{listing}{listings}
\Crefname{sublstlisting}{Listing}{Listings}

\hypersetup{
  colorlinks   = true, %Colours links instead of ugly boxes
  urlcolor     = blue, %Colour for external hyperlinks
  linkcolor    = blue, %Colour of internal links
  citecolor   = red %Colour of citations
}

\pgfplotsset{compat=1.17}
\def\tablename{Table}

\usepackage{tikz}
\usetikzlibrary{arrows.meta}
\usepackage{booktabs}   %% For formal tables:
\usepackage{subcaption} %% For complex figures with subfigures/subcaptions
\usepackage{comment}

\usepackage{amsfonts}

\usepackage{xcolor, color, soul}
\usepackage{adjustbox}
\sethlcolor{green}

\newcounter{sema}
\DeclareRobustCommand{\rsema}[1]{%
   \refstepcounter{sema}S\thesema\label{#1}%
   }

\newcommand{\cladtool}{\textsc{CHEF-FP}\xspace}
\RequirePackage[normalem]{ulem} %DIF PREAMBLE
\RequirePackage{color}\definecolor{RED}{rgb}{1,0,0}\definecolor{BLUE}{rgb}{0,0,1}\definecolor{purple}{rgb}{0.5,0,0.5} %DIF PREAMBLE

\DeclareRobustCommand*{\IEEEauthorrefmark}[1]{\raisebox{0pt}[0pt][0pt]{\textsuperscript{\footnotesize\ensuremath{
  \ifcase#1
  \relax\or
  \star\or
  \dagger\or
  \ddagger\or
  \mathsection
  \else\textsuperscript{\expandafter\romannumeral#1}
  \fi
}}}}

\newcommand\copyrighttext{%
  \footnotesize \textcopyright 2012 IEEE. Personal use of this material is permitted.
  Permission from IEEE must be obtained for all other uses, in any current or future
  media, including reprinting/republishing this material for advertising or promotional
  purposes, creating new collective works, for resale or redistribution to servers or
  lists, or reuse of any copyrighted component of this work in other works.}
\newcommand\copyrightnotice{%
\begin{tikzpicture}[remember picture,overlay]
\node[anchor=south,yshift=10pt] at (current page.south) {\fbox{\parbox{\dimexpr\textwidth-\fboxsep-\fboxrule\relax}{\copyrighttext}}};
\end{tikzpicture}%
}

\begin{document}

 \author{ 
  \IEEEauthorblockN{
  Garima Singh\IEEEauthorrefmark{1}\IEEEauthorrefmark{4}, 
  Baidyanath Kundu\IEEEauthorrefmark{1}\IEEEauthorrefmark{4}, 
  Harshitha Menon\IEEEauthorrefmark{2}, 
  Alexander Penev\IEEEauthorrefmark{3},
  David J. Lange\IEEEauthorrefmark{4},
  Vassil Vassilev\IEEEauthorrefmark{4}\IEEEauthorrefmark{1}\\~\\}

  \IEEEauthorblockA{\itshape
  \IEEEauthorrefmark{1}European Council for Nuclear Research, Espl. des Particules 1, 1211 Meyrin, Switzerland \\
  \IEEEauthorrefmark{4}Department of Physics, Princeton University, Princeton, New Jersey 08544, USA \\
  \IEEEauthorrefmark{2}Center for Applied Scientific Computing, Lawrence Livermore National Laboratory, Livermore, California 94551 USA \\
  \IEEEauthorrefmark{3}Faculty of Mathematics and Informatics, University of Plovdiv, 236 Bulgaria Blvd., 4000 Plovdiv, BULGARIA \\
  E-mail: \IEEEauthorrefmark{1}\{garima.singh, baidyanath.kundu, vassil.vassilev\}@cern.ch, \IEEEauthorrefmark{2}harshitha@llnl.gov, \IEEEauthorrefmark{3}apenev@uni-plovdiv.bg, \\
  \IEEEauthorrefmark{4}david.lange@princeton.edu}
}

%%%\begin{document}
% \input{coverLetter}
% \newpage
%% Title information
\title{Fast And Automatic Floating Point Error Analysis With \cladtool}
%\thanks{Submitted to the editors \today.\\ This work was performed under the auspices of the U.S. Department of Energy by Lawrence Livermore National Laboratory under Contract DE-AC52-07NA27344. Work at LLNL was funded by the Laboratory Directed Research and Development Program under project tracking code 20-ERD-043, LLNL-JRNL-XXXXXX-DRAFT.}}
%\author

%%%\makeatletter
%%%\newcommand{\linebreakand}{%
%%%  \end{@IEEEauthorhalign}
%%%  \hfill\mbox{}\par
%%%  \mbox{}\hfill\begin{@IEEEauthorhalign}
%%%}
%%%\makeatother

%\author{
% \IEEEauthorblockN{Garima Singh}
% \IEEEauthorblockA{\textit{CERN}}
% \IEEEauthorblockA{Princeton University}
% 
% \and
%
% \IEEEauthorblockN{Baidyanath Kundu}
% \IEEEauthorblockA{\textit{CERN}}
% \IEEEauthorblockA{Princeton University}
%
% \and
%
% \IEEEauthorblockN{Harshitha Menon}
% \IEEEauthorblockA{\textit{Lawrence Livermore National Laboratory}}
%
% %\and
% %\newline
% \linebreakand
%
% \IEEEauthorblockN{Alexander Penev}
% \IEEEauthorblockA{\textit{ University of Plovdiv Paisii Hilendarski}}
%
% \and
%
% %\newline
% \IEEEauthorblockN{David J. Lange}
% \IEEEauthorblockA{Princeton University}
%
% \and
%
% \IEEEauthorblockN{Vassil Vassilev}
% \IEEEauthorblockA{Princeton University}
% % }

\thispagestyle{plain}
\pagestyle{plain}

\maketitle

\copyrightnotice

% \mynote{Harshitha: Write an abstract to submit for the 29th.}
\begin{abstract}
As we reach the limit of Moore’s Law, researchers are exploring different paradigms to achieve unprecedented performance. Approximate Computing (AC), which relies on the ability of applications to tolerate some error in the results to trade-off accuracy for performance, has shown significant promise. Despite the success of AC in domains such as Machine Learning, its acceptance in High-Performance Computing (HPC) is limited due to its stringent requirement of accuracy.
%in HPC applications, blindly applying AC techniques throughout the code is not a viable option as it would result in output errors beyond a prescribed threshold. 
We need tools and techniques to identify regions of the code that are amenable to approximations and their impact on the application output quality so as to guide developers to employ selective approximation. To this end, we propose \cladtool{}, a flexible, scalable, and easy-to-use source-code transformation tool based on Automatic Differentiation (AD) for analysing approximation errors in HPC applications.

\cladtool uses Clad, an efficient AD tool built as a plugin to the Clang compiler and based on the LLVM compiler infrastructure, as a backend and utilizes its AD abilities to evaluate approximation errors in C++ code. \cladtool{} works at the source level by injecting error estimation code
into the generated adjoints. This enables the error-estimation code to undergo compiler optimizations resulting in improved analysis time and reduced memory usage. 
We also provide theoretical and architectural augmentations to source code transformation-based AD tools to perform FP error analysis.
In this paper, we primarily focus on analyzing errors introduced by mixed-precision AC techniques, the most popular approximate technique in HPC. We also show the applicability of our tool in estimating other kinds of errors by evaluating our tool on codes that use approximate functions. Moreover, we  demonstrate the speedups achieved by CHEF-FP during analysis time as compared to the existing state-of-the-art tool as a result of its ability to generate and insert approximation error estimate code directly into the derivative source. The generated code also becomes a candidate for better compiler optimizations contributing to lesser runtime performance overhead.

% Computer architecture supports multiple levels of floating-point precision, and mixed precision computing is the use of different numerical precision in a program in order to improve the performance (runtime and/or energy) of an application. One of the significant challenges to applying mixed precision in HPC is identifying where to apply mixed precision and understanding how it will impact the application output quality. Our work uses Clad, an efficient automatic-differentiation tool built as a plugin to the Clang compiler and based on the LLVM compiler infrastructure. We use Clad as a backend for the presented framework and utilize its AD abilities to evaluate floating-point errors in C++ code. Finally, because we generate and inject FP error estimation code directly into the derivative source, we notice a significant speedup compared to other tools performing similar FP error analysis. The generated code also becomes a candidate for better compiler optimizations contributing to lesser runtime performance overhead.
\end{abstract}

%% Keywords
%% comma separated list
%\keywords{ Approximate computing, high-performance computing, Clang, Automatic Differentiation}  %% \keywords are mandatory in final camera-ready submission
%%%\begin{IEEEkeywords}Approximate computing, high-performance computing, Clang, Clad, Automatic Differentiation \end{IEEEkeywords}

%% \maketitle
%% Note: \maketitle command must come after title commands, author
%% commands, abstract environment, Computing Classification System
%% environment and commands, and keywords command.
%%%\maketitle

% REQUIRED
%\begin{AMS}
%  65G50, % Roundoff error
%  65Y04 % Numerical algorithms for computer arithmetic
%\end{AMS}

% Body
\section {Introduction}
\label{sec:introduction}
%\mynote{Garima: Profile why we take less memory and if that's the case put that as a statement here.}
As we enter the post-Moore era, where we no longer enjoy the free lunch of performance growth from shrinking the transistor features, researchers are exploring other computing paradigms to increase computational throughput. Approximate Computing (AC) has garnered significant interest as a promising approach for increasing peak performance. AC relies on the application to tolerate some amount of error to achieve performance gains. Among the various AC techniques that currently exist, reduced floating-point (FP) precision, or mixed precision, has gained in popularity.
Computer architectures support multiple levels of precision for FP data and arithmetic operations — 64 bits \emph{double} precision,
32 bits \emph{single} precision, 128 bits \emph{quad} precision, and 16 bits \emph{half} precision. The choice
of precision determines the amount of rounding error. While using higher precision for data and operations may result in increased accuracy, it can lead
to an increase in application execution time, memory and energy consumption. Mixed-precision tuning involves
using higher precision when necessary to maintain accuracy and using lower precision where
we can improve performance.

Despite the availability of multiple levels of precision for FP, it is 
challenging to apply them in HPC and scientific codes. To use them effectively,
developers need to understand the details of rounding errors as well as 
how they propagate through their applications.
Due to the lack of scalable and rigorous tools that analyze
error sensitivity in the applications’ code regions, developers often
resort to the safer option of using high precision 
%(\emph{double})
throughout. 
We need scalable tools to understand the impact of lowering the precision of 
data and computation, identify error-tolerant regions, and provide guidance to users.

% Floating-point arithmetic operations are the most prevalent computation
% in HPC applications. Computer architectures support multiple levels of
% precision for floating-point (FP) data and arithmetic operations. In the IEEE 754 standard, the
% most common representation today for real numbers, choices are 64 bits \emph{double} precision,
% 32 bits \emph{single} precision, 128 bits \emph{quad} precision, and 16 bits \emph{half} precision. The choice
% of precision determines the amount of rounding error. 
% Therefore, we need to be prudent about how and where lower-precision computing is
% applied in a HPC code.

Several techniques have been proposed to estimate the sensitivity profile of an application,
including automated search based approaches~\cite{rubio2013precimonious,lam2018fine}, static analysis~\cite{das2021robustness, solovyev2018rigorous}, or using Automatic Differentiation (AD)~\cite{menon2018adapt}.
Unfortunately, the existing set of tools fail to provide a feasible solution for HPC applications. 
Search-based approaches are very expensive as the
state space is significantly large and quickly become infeasible even for small benchmarks.
Static analysis-based approaches using
interval analysis or Taylor series approximation provide rigorous estimates for FP
errors but are so far limited to programs with a small number of operations~\cite{solovyev2018rigorous}. 
Several methods~\cite{menon2018adapt, vassiliadis2016towards} leverage AD to determine error-resilient regions of codes.
While they have been shown to work well for smaller HPC benchmarks, they often require
manual code changes and incorporating several software tools together~\cite{lam2019floatsmith}.  
These tools are slow and have high memory overhead, making them infeasible for large HPC workloads.

% Stuff added by Garima 

% Our contribution builds on the work by Harshitha Menon et al. on ADAPT. ADAPT is a floating-point (FP) error estimation tool that efficiently delivers mixed-precision configurations for C++ codes using operator overloading based automatic differentiation (AD). While ADAPT's approach is succinct and practical for smaller workloads,  it exposes a set of points for optimization and automatization on larger projects. 
We propose \cladtool{}, a scalable, flexible, and easy-to-use tool, for analyzing approximation errors in HPC applications.
We use Clad~\cite{clad}, an efficient AD tool built as a plugin to the Clang compiler
%and based on the LLVM compiler infrastructure, 
as a backend for the presented framework and utilize its AD abilities to evaluate FP errors in C++ code. %Clad implements source code transformation based AD, i.e., it inspects the internal compiler representation of a program to generate its derivative.
Clad builds derivatives by transforming the internal compiler representation of the program.
This level of granularity allows \cladtool{}, built as an extension to Clad, to exploit a program's source information to automatically add auxiliary instructions to the generated derivative depending on a set of rules. Accordingly, \cladtool{} annotates code with FP error information without user intervention. This automation is beneficial for large-scale codes wherein going into hand-write annotations becomes complex, tedious, and error-prone.
Moreover, \cladtool{} enables domain-specific FP error analysis due to the increased flexibility of a compiler-based backend. This allows easy modification of parts of the \cladtool{} framework to achieve a more tailored FP error analysis. This level of variability 
%also
enables a  finer-grained analysis where users can themselves query source information and make more precise decisions for the direction of the analysis.     

Finally, because we generate FP error estimation (EE) code directly into the derivative source, we observe a significant speedup compared to other tools performing similar FP error analysis. The generated code also becomes a candidate for further compiler optimizations contributing to better runtime performance.
We use the FP error profile provided by \cladtool{} to guide the design of a mixed-precision version of the code and achieve performance improvements of $8\%$ for \textit{HPCCG} and $65\%$ for \textit{Black-Scholes} while satisfying the user-specified error threshold. At analysis time, \cladtool obtained a maximum time speedup of 217\% with the Simpsons benchmark and a memory efficiency of 632\% with the Black-Scholes benchmark.

% \mynote{Compare with existing works and highlight their shortcomings.}

% \mynote{Add a list of contributions.}
\textbf{Key Contributions}
%\textcolor{cyan}{The contribution of our work is a theoretical and architectural enhancement of a compiler-based, source-transformation approach to AD in an effort to aid FP EE in a loosely-coupled way such that it does not overwhelm the already overly complex reverse accumulation AD mode. Traditionally, AD intertwines analyses reducing the tracked variables hindering AD-based FP EE. The theoretical contribution extends the widely-adopted AD operational semantics notation. The implementation extends Clad by adding support for third-party plugins. We implemented CHEF-FP as a loosely-coupled listener to AD generation events. Naively, CHEF-FP "hijacks" the gradient generation process to augment the adjoints with FP-related code. Clad provides the AD infrastructure for the gradients, CHEF-FP generates code for caching values, accumulating errors, tracking derivatives, and other information directly into the gradients. CHEF-FP captures structural and source information allowing better analysis of loops and nested calls. User-wise, CHEF-FP offers a simple interface allowing users to express analysis code in C++, enabling straightforward analyses to avoid writing complex compiler/AD API calls. In particular we show:}
\begin{itemize}
    \item \cladtool{}, an efficient tool to automate AD-based FP error analysis. CHEF-FP inlines error calculations into the adjoint code and results in faster and more memory-efficient analysis, making it suitable for use in data-intensive applications;
    %The output of \cladtool{} guides developers to design a mixed-precision version of code that satisfies a given error threshold;
    \item Formalism for augmenting AD to perform FP error analysis. Specifically, we present a generic way to extend source code transformation-based AD tools;
    \item Customizable EE module that supports any AD-based user-defined error model;
    \item Tool evaluation using a set of HPC benchmarks for mixed-precision and approximate function error analysis.
\end{itemize}
%\section{Overview}
%\mynote{Provide an overview of the approach with an example highlighting 
%the shortcomings of the previous approach. Talk about the current flow (step1, step2, step3) and demonstrate what we can do in a single step. Eg, explain the problem and tell we have a solution.}

\section{Background}
Floating-point arithmetic operations are the most prevalent computation
in HPC applications. Computer architectures support multiple levels of
precision for FP data and arithmetic operations. In the IEEE 754 standard \cite{8766229}, the
most common representation today for real numbers, choices are 64 bits \emph{double} precision,
32 bits \emph{single} precision, 128 bits \emph{quad} precision, and 16 bits \emph{half} precision. The choice
of precision determines the amount of rounding error.
% Most real-valued applications heavily rely on FP arithmetic. However, due to the finite nature of computing, it is incredibly hard to exactly represent a lot of the real values as FP numbers. This inexactness inevitably leads to errors that can grow quickly and have serious implications on the accuracy of the application. These errors born out of squeezing the infinite real-valued numbers into a finite set of bits are known as \emph{rounding errors}.

The FP rounding error is the accumulation of FP errors from each variable toward the target function's result. Target functions are modeled as a series of multiple assignment operations. These assignments and function inputs are assumed to be independent (FP error analysis on correlated variables is still a nascent field, and definitive methods of estimation are yet to be established). We define an expression to calculate and accumulate the FP error introduced by each assignment and refer to these error calculation and accumulation expressions as FP error models. \cladtool{}, %with its default error model, 
aims to estimate an upper bound on these FP rounding errors. 

\subsection{Modelling Floating-Point Errors}
\label{sec:fpmodel}

The error model describes a metric that can be used in error estimation analysis. This metric is applied to all assignment operations in the function, and the results from its evaluation are accumulated into the total FP error of the left-hand side of the assignment. The Taylor series approximation is well suited to model floating point errors. Let's assume an arbitrary function \(y = f(x)\), where \(x\) is represented in the standard IEEE 754 single precision. Assuming a floating point error of \(h\) in \(x\), we define \(\widetilde{f}(x) = f(x + h) \). A symbolic Taylor series expansion yields:% the following:
\[ f(x + h) = f(x) + \frac{h}{1!}f'(x) + \frac{h^2}{2!}f''(x) + \frac{h^3}{3!}f'''(x) + ...
\]

Here \(f'(x)\) represents the derivative of the function with respect to \(x\). Approximating the series to
%its
first-order results in: % the following:
\[ f(x + h) = f(x) + \frac{h}{1!}f'(x) + O(h^2)
\]

An expression for the absolute floating point error in \(f\) is:
%can be represented as:
\[ A_f = | \widetilde{f}(x) - f(x) |
   => A_f = | \frac{h}{1!}f'(x) | 
   => A_f = | hf'(x) |
\]

To determine the maximum absolute floating point error in \(f\), we write \(h = \epsilon_m |x|\) where \(\epsilon_m\) is known as the \textit{machine epsilon}. The machine epsilon gives the maximum relative representation error in floating point variables due to rounding. It is a machine-dependent value and follows the IEEE standard in most compilers. The absolute error can then be written as:
\begin{equation}
 \label{eqn:deftmodel}    
 A_f=|\epsilon_m|x|f'(x)|
\end{equation}
 
This error model is sufficient for most smaller cases and can produce loose upper bounds of the maximum permissible FP error in programs. %\mynote{Expand.}
%FIXME: Move to implementation details;We provide a simple class interface that users can implement to develop their error models. We also provide comprehensive demos on how to code up a custom error model from scratch.
Here, AD techniques can enable efficient computations of $f'(x)$ that scale to real-world workflows.

\subsection{Automatic Differentiation Basics}

%Automatic differentiation uses the properties of the chain rule to evaluate differential operators of a function that is represented as a computer program. A differential operator maps an input function to a new function that captures how small perturbations in the inputs affect the output. First-order differential operators of real-valued functions include directional derivatives and gradients~\cite{betancourt2018geometric}. A differential operator, or a \textit{pushforward}, maps between the local tangent space of the input (the space of perturbations to the inputs) to the local tangent space of the outputs. Informally, a pushforward takes a wiggle in the input space, and tells what wobble is created in the output space, by passing it through the function. The \textit{pullback} operator maps the cotangent space of the outputs to the inputs. Informally, the pullback takes wiggles with respect to the function’s output, and tells the equivalent wobble with respect to the functions input.

AD takes as input program code that \textit{has meaningful differentiable properties} and produces new code augmented with pushforward or pullback operators~\cite{betancourt2018geometric}. The AD mode, which produces pushforward operators to capture the sensitivity from inputs to outputs, is commonly known as \textit{forward} mode. The AD mode capturing the sensitivity from outputs to inputs is called \textit{reverse} mode, \textit{adjoint} mode, or backpropagation. The computational merit of the pullback is that it provides a very eﬀicient way to compute the function's gradient with relative time complexity which is independent of the input size~\cite{griewank2008evaluating}.

AD tools can be categorized by how much work is done before program execution. AD typically starts by \textit{building a computational graph}, or a directed acyclic graph of mathematical operations applied to an input. At one extreme, the \textit{tracing}, or \textit{taping} approach constructs and processes the computational graph at the time of execution each time a function is invoked. In contrast, the source transformation approach does as much as possible at compile time to create the derivative only once. The tracing approach is easier to implement and adopt into existing codebases. The source transformation approach has better performance but usually covers a subset of the language.
\textbf{Tracing:} Records the linear sequence of
computation operations at runtime into a \textit{tape}. The control flow is flattened to produce a derivative. A typical implementation is via operator overloading, defining a special floating type with overloaded elementary  operations. Algorithms use this type to trigger   differentiation by calling a special function. There are numerous C++ AD tools  based on tracing, including ADOL-C~\cite{griewank1996algorithm}, and Adept~\cite{hogan2014fast}. As the derivative is produced at runtime, the just-in-time differentiation process is constrained to perform optimizations quickly.
Metaprogramming techniques, such as expression templates, can mitigate the issue, but they cannot optimize across  statements and generally do not handle control flow~\cite{hogan2014fast, sagebaum2018expression}.
\textbf{Source Transformation:} Constructs  the computation graph and produces a derivative at compile time. More compile-time optimizations can be applied, such as reorganizing or evaluating simple constant expressions and common subexpression elimination. Source transformation is more difficult to implement as it requires a significant investment in developing and maintaining a language parser. Tapenade~\cite{hascoet2013tapenade} is an example of a source transformation tool with custom parsers for C and Fortran. Source transformation tools usually do not support the full language feature set (e.g. certain language idioms are particularly hard to differentiate).

Historically, toolmakers made trade-offs between ease of use, performance, and ease of integration. AD now benefits from better language support to avoid such trade-offs. Recently, production compilers like Clang allowed tools to reuse the language parsing infrastructure. Enzyme~\cite{enzymeNeurips} and Clad~\cite{clad} are examples of compiler-based AD tools using such preexisting parsers.

\section{AD-based FP Error estimation Using CHEF-FP}

%The reverse accumulation mode transforms the original function into two generic steps called a \emph{forward sweep} and a {backward sweep}. During the forward sweep the original function augmented with statements keep track of the state changes allowing the function to be executed backwards in its backward sweep. The backward sweep is responsible for computing the partial derivatives or \emph{adjoints}.

%In turn, the EE framework augments the backward sweep by instrumenting every assignment which is used in an adjoint, giving control to a special function which implements the specific EE model. Then in the end 

An important aspect of dealing with FP applications with high precision requirements is identifying the sensitive, or more error-prone, areas to devise suitable mitigation strategies. This is where AD-based sensitivity analysis can be very useful. It can find specific variables or regions of code with a high contribution to the overall FP error in the application. 

The sensitivity of a variable \(x\) to FP errors (\(S_x\)) can be deduced from the default error model described in Eq. \ref{eqn:deftmodel} as:
\[S_x = | x f'(x) |\]

Here, the total contribution of FP errors to the function (either a routine or an entire program) by variable \(x\) increases as \(S_x\) increases. A study of the trends of these sensitivity values across the domain of a function can reveal insights into the numerical stability of the function and help determine possible causes of instabilities. Sensitivity values can also be used as a guide for a class of type-based optimizations called \emph{Mixed Precision Tuning}.

Mixed Precision Tuning involves \textit{demoting} certain variables to lower precision without severely affecting the overall accuracy of the application. In corollary, it involves preserving or \textit{promoting} the precision of variables that have a significant effect on the accuracy. A mixed precision tuned configuration is only valid when the difference of the pre- and post-tuning accuracy is less than some defined threshold value. An effective way to maintain this requirement is by analyzing the sensitivity of all input and intermediate variables and selecting the ones with lower sensitivity to be demoted. The FP error contributions of the demoted variables are accumulated and compared to the threshold value. A mixed precision configuration is reached when the accumulated error meets the threshold value.

Taylor-based analyses described in~\cref{sec:fpmodel} require modeling assignments to FP variables and their respective adjoints. The adjoint accumulation mode naturally offers such mapping. The \cladtool implementation exports this information from Clad. Clad can be used either as a part of the compilation lowering pipeline or to generate source code that can be compiled by another compiler toolchain. Clad implements forward and adjoint mode AD, together with a flexible extension system that allows user code to subscribe to events during the process of adjoint creation.

Compiler-based AD tools can operate at the level of different program representations, and each implementation has its own set of pros and cons. For example, Clad implements AD on Clang's high-level representation to make use of better diagnostics, support compile-time programming, generate understandable source code, and use a standard optimization pipeline. These properties are key for AD-based FP analyses. If AD runs before the optimization pipeline, unsafe optimizations might introduce extra floating point errors. Running AD after optimization, as done by tools including Enzyme, avoids these errors. In this case, even standard optimizations could break differentiability. It is still to be seen if there is a way to combine the strengths of both approaches.

While many AD-based approaches exist for error analysis, they typically involve significant manual effort to perform code changes or long toolchains to automate it. In this work, we present \textit{\cladtool}, an AD-based FP error estimation framework that requires less manual integration work and comes packaged with Clad, taking away the tedious task of setting up long toolchains. The tool's proximity to the compiler and the fact that the FP error annotations are built into the code's derivatives allows for powerful compiler optimizations that provide significant speedups of the analysis time when compared to the current state-of-the-art tools like ADAPT~\cite{menon2018adapt}.

Another advantage of using source-level AD tools (like Clad) as the backend is that they provide important source insights. A higher-level representation, for example, enables us to identify and attribute FP errors to variables as they are visible at the source. Recovery of such information (variable names, IDs, source location, and so on) becomes difficult for lower-level tools, where optimizations can optimize away variables of interest. Clad can also identify special constructs (such as loops, lambdas, functions, if statements, and so on) and tune the error code generation accordingly. This information, while available to lower-level tools, is more difficult to extract and process at that level.

\subsection{Source Transformation FP Error Estimation Using AD}

Algorithm~\ref{lst:algo} describes a transformation that connects the variable-adjoint mapping available in a source transformation AD engine with an FP EE framework. For each function annotated for EE, we use the mapping (in \textit{AdjointAD}) while differentiating (\textit{NewFunction}) to hand the control to the EE model (\textit{AssignError}) to insert extra instructions. 
\begin{algorithm}
	\caption{Error Estimation Generation Process} 
	\begin{algorithmic}[1]
	    \Require Selected by $estimate\_error$ functions
	    \Ensure Generated error estimated functions "$\overline{function}$"
		\ForAll {$function$ in $\textit{estimate\_error.functions}$}
		    \State $map\langle variable, adjoint \rangle$ $\gets$ \Call{AdjointAD}{$function$}
		    \State $\overline{function}$ $\gets$ \Call{NewFunction}{$function, map$}
			\ForAll {$variable,adjoint$ in $map$}
				\State \Call{AssignError}{$\overline{function}, variable, adjoint$}
%				    \Comment{This is a comment}
			\EndFor
			\State \Call{FinalizeEE}{$\overline{function}$}
		\EndFor
	\end{algorithmic}
 \label{lst:algo}
\end{algorithm}
In the end, we compute the total error (\textit{FinalizeEE}) by allowing the custom model to pass it as an output parameter or print it on the screen. One important aspect is that the algorithm does not impose constraints on the structure of the right-hand side of the FP assignments. For implementations of certain functions, this implies that the right-hand side can contain arbitrary long expressions. In turn, the error contribution of longer expressions is computed with less precision. This feature gives users implicit control by allowing them to reduce the expression size manually. The generated code depends on the particular implementation of the EE  model; however, we can schematically list a possible generalization for a better understanding of the process.

\subsection{Sample Syntactic Structure of the Generated Function}

Algorithm~\ref{lst:algo} transforms the input program in Fig~\ref{fig:CladOriginalFnStructure} and produces another program as shown in Fig~\ref{fig:CladEEStructure}. 
Fig~\ref{fig:CladOriginalFnStructure} represents a function written in a computer programming language that takes $P$ parameters and returns a result of type $T$. Without loss of generality, its body consists of instructions denoted by $L_i,\forall i\in[1..n]$, where $S=f_i(S)$ manipulating internal state $S$ and returning the final state $S$ in $result(S)$. This notation does not exclude control flow constructs which can be represented with a sufficiently long linear sequence of $L_i$.

The adjoint accumulation mode of AD computes partial derivatives starting from the function's ($\overline{FuncName}$) outputs 
towards the function inputs ($P$). It requires executing the function in reverse order.
Fig.~\ref{fig:CladEEStructure} illustrates the AD adjoint accumulation transformation, which is a possibly augmented \begin{figure}[h]
    \centering
    {\renewcommand{\arraystretch}{1.2}
    \begin{tabular}{ |p{1.75em} l| }
      \arrayrulecolor{black} \hline
      \multicolumn{2}{|l|}{$\textcolor{Blue}{\texttt{function}}\;FuncName\textcolor{Blue}{(}P\textcolor{Blue}{):}\;T$} \\
      
      \arrayrulecolor{black} \hline
      
      {} & $Initialize\;S\;with\;P$ \\
      \arrayrulecolor{lightgray} \hline
      $L_1$: & $S\;\textcolor{Blue}{=}\;f_1\textcolor{Blue}{(}S\textcolor{Blue}{)}$ \\
      {} & $\ldots$ \\
      $L_i$: & $S\;\textcolor{Blue}{=}\;f_i\textcolor{Blue}{(}S\textcolor{Blue}{)}$ \\
      {} & $\ldots$ \\
      $L_{n}$: & $S\;\textcolor{Blue}{=}\;f_{n}\textcolor{Blue}{(}S\textcolor{Blue}{)}$ \\
      
      \arrayrulecolor{lightgray} \hline
      
      {} & $\textcolor{Blue}{\texttt{return}}\;result(S)$ \\

      \arrayrulecolor{black} \hline
    \end{tabular}
    }

    \caption{Structure of the original function $FuncName$}
    \label{fig:CladOriginalFnStructure}
    \vspace{-0.0cm}
\end{figure}
instruction sequence $L_i$. The augmented instruction records a subset of the internal state (denoted as $out(L_i) \subset S$, which depends on $f_i$) necessary to preserve the semantics
\begin{figure}[h]
    \centering
    {\renewcommand{\arraystretch}{1.2}
    \begin{tabular}{ |p{1.75em} l!{\color{lightgray}\vline}c| }
      \arrayrulecolor{black} \hline
      \multicolumn{3}{|l|}{$\textcolor{Blue}{\texttt{function}}\;\overline{FuncName}\textcolor{Blue}{(}P\textcolor{Blue}{,}\;E\textcolor{Blue}{):}\;void$} \\
      
      \arrayrulecolor{black} \hline
      
      {} & \multicolumn{2}{l|}{$Initialize\;S\;and\;\overline{S}$} \\
      \arrayrulecolor{lightgray} \hline
      $L_1$: & $\textcolor{OliveGreen}{\texttt{Push}}\textcolor{Blue}{(}out\textcolor{Blue}{(}L_1\textcolor{Blue}{));}\;S\;\textcolor{Blue}{=}\;f_1\textcolor{Blue}{(}S\textcolor{Blue}{)}$ & \multirow{5}{*}{\rotatebox[origin=c]{90}{\small{\emph{\textcolor{gray}{forward sweep}}}}} \\
      {} & $\ldots$ & \\
      $L_i$: & $\textcolor{OliveGreen}{\texttt{Push}}\textcolor{Blue}{(}out\textcolor{Blue}{(}L_i\textcolor{Blue}{));}\;S\;\textcolor{Blue}{=}\;f_i\textcolor{Blue}{(}S\textcolor{Blue}{)}$ & \\
      {} & $\ldots$ & \\
      $L_{n-1}$: & $\textcolor{OliveGreen}{\texttt{Push}}\textcolor{Blue}{(}out\textcolor{Blue}{(}L_{n-1}\textcolor{Blue}{));}\;S\;\textcolor{Blue}{=}\;f_{n-1}\textcolor{Blue}{(}S\textcolor{Blue}{)}$ & \\
      
      \arrayrulecolor{lightgray} \hline
      
      $\overleftarrow{L_n}$: & $\overline{S}\;\textcolor{Blue}{=}\;\overline{S} \textcolor{Blue}{\times} f'_n\textcolor{Blue}{(}S\textcolor{Blue}{);}\; \textcolor{OliveGreen}{\texttt{AssignError}}\textcolor{Blue}{(}\overleftarrow{L_n}\textcolor{Blue}{)}$ & \multirow{8}{*}{\rotatebox[origin=c]{90}{\small{\emph{\textcolor{gray}{backward sweep}}}}} \\
      {} & $\textcolor{OliveGreen}{\texttt{Pop}}\textcolor{Blue}{(}out\textcolor{Blue}{(}L_{n-1}\textcolor{Blue}{))}$ & \\
      {} & $\ldots$ & \\
      {} &  $\textcolor{OliveGreen}{\texttt{Pop}}\textcolor{Blue}{(}out\textcolor{Blue}{(}L_i\textcolor{Blue}{))}$ & \\
      $\overleftarrow{L_i}$: & $\overline{S}\;\textcolor{Blue}{=}\;\overline{S} \textcolor{Blue}{\times} f'_i\textcolor{Blue}{(}S\textcolor{Blue}{);}\; \textcolor{OliveGreen}{\texttt{AssignError}}\textcolor{Blue}{(}\overleftarrow{L_i}\textcolor{Blue}{)}$ & \\
      {} & $\ldots$ & \\
      {} & $\textcolor{OliveGreen}{\texttt{Pop}}\textcolor{Blue}{(}out\textcolor{Blue}{(}L_1\textcolor{Blue}{))}$ & \\
      $\overleftarrow{L_1}$: & $\overline{S}\;\textcolor{Blue}{=}\;\overline{S} \textcolor{Blue}{\times} f'_1\textcolor{Blue}{(}S\textcolor{Blue}{);}\; \textcolor{OliveGreen}{\texttt{AssignError}}\textcolor{Blue}{(}\overleftarrow{L_1}\textcolor{Blue}{)}$ & \\
      
      \arrayrulecolor{lightgray} \hline \arrayrulecolor{black}

      {} & \multicolumn{2}{l|}{$E\;\textcolor{Blue}{=}\;\textcolor{OliveGreen}{\texttt{FinalizeEE}}\textcolor{Blue}{(}\overline{FuncName}\textcolor{Blue}{,}\;\overline{S}\textcolor{Blue}{)}$} \\

      \arrayrulecolor{lightgray} \hline \arrayrulecolor{black}
%      {} & \multicolumn{2}{l|}{$\textcolor{Blue}{\texttt{return}}\;result(S)$} \\

      \arrayrulecolor{black} \hline
    \end{tabular}
    }

    \caption{Structure of the error estimated function $\overline{FuncName}$}
    \label{fig:CladEEStructure}
\end{figure}
when the instructions are executed in reverse order. A common implementation mechanism is to use a LIFO structure such as a \emph{stack} and insert \emph{push} or \emph{pop} operations via $\textit{Push}$ or $\textit{Pop}$. This transformation is known as \emph{forward sweep}. The \emph{backward sweep} is responsible for computing the adjoint ($\overline{S}$) for each instruction $\overleftarrow{L_i}$. $\overline{S}$ may require restoring altered state $S$ by using a \emph{pop} operation to correctly evaluate $f_i'(S)$.

Our adjoint accumulation mode extension adds 3 elements: another output parameter $E$ to $\overline{FuncName}$ modeling the total error; a callback call to $\textit{AssignError}$ for every assignment taking variable and its adjoint as parameters; and a function $\textit{FinalizeEE}$ computing the total error $E$.

\subsection{Structural Operational Semantics of the Transformation}

%Fig.~\ref{fig:CladOriginalFnStructure} and Fig~\ref{fig:CladEEStructure} can be described similarly to the structural operational semantics notation used by Hascoet and Pascual in~\cite{hascoet2013tapenade}. In brief, the notation represents a logical implication where the rules on top of the fraction is the \emph{antecedent} below is the \emph{consequent}. Due to space limitations we show only essential rules from the FP error estimation perspective. 
The transformation from the pseudo-code shown in Fig.~\ref{fig:CladOriginalFnStructure} to  its error estimation form, shown in Fig~\ref{fig:CladEEStructure}, can be described with the structural operational semantics notation. Similar to Hascoet and Pascual's work describing AD semantics in~\cite{hascoet2013tapenade}, we can extend the rules to describe the EE specific transformations. Due to space limitations, we  show only essential rules supporting FP EE. They describe how a generic source-transformation AD framework can be extended to support AD-based FP analysis. The rest of the rules associated with AD are already available in~\cite{hascoet2013tapenade}.
\[
\inference[\texttt{\rsema{sema:FuncDecl}:}]
%\inference[(1)]
{
  \textcolor{OliveGreen}{\textit{isEstErrFunction}}(FuncName) \\
   FuncName \xrightarrow {newFunction} \overline{FuncName} \\
  Params -> \overline{Params} \quad
  Locals, Params -> \overline{Locals} \\
  Stmts ->
    \left[
      \makecell[l]{
        \overrightarrow{Stmts} \\
        \overleftarrow{Stmts}
      }
    \right. \\
}
{
  \makecell[l]{
    \textcolor{Blue}{\texttt{function}}\;FuncName\textcolor{Blue}{(}Params\textcolor{Blue}{)}\;\textcolor{Blue}{:}\;T \; \textcolor{Blue}{\{} \\
    \qquad Locals\textcolor{Blue}{;}\;Stmts \\\textcolor{Blue}{\}} \\ 
      -> \textcolor{Blue}{\texttt{function}}\;\overline{FuncName}\textcolor{Blue}{(}\overline{Params}\textcolor{Blue}{,}\;E\;\textcolor{Blue}{:}\;T\textcolor{Blue}{)}\;\textcolor{Blue}{:}\;void \; \textcolor{Blue}{\{} \\
       \qquad \quad \overline{Locals}\textcolor{Blue}{;}\; \overrightarrow{Stmts}\textcolor{Blue}{;}\;\overleftarrow{Stmts}\textcolor{Blue}{;}\; \\
       \qquad \quad E\;\textcolor{Blue}{\texttt{=}}\; \textcolor{OliveGreen}{\texttt{FinalizeEE}} \textcolor{Blue}{(}FuncName\textcolor{Blue}{)} \\
       \quad \; \textcolor{Blue}{\}}
  }
}
\]

In rule S1, %TODO: Use ref S\ref{sema:FuncDecl}
the premises are fulfilled when:
\begin{enumerate*}[label=(\roman*)]
    \item Predicate $\textit{isEstErrFunction}$ is true when \textit{FuncName} is selected for being an error estimation candidate;
    \item \textit{newFunction} can create a $\overline{FuncName}$;
    \item \textit{Params} and \textit{Locals} can be transformed; and
    \item The original function's statement sequence has been transformed into a forward and backward sweep.
\end{enumerate*}
When all antecedents are fulfilled, $\textit{FinalizeEE}$ is added at the end of the backward sweep. 
%When the predicate $\textit{isEstErrFunction}$ is true and the original function's statement sequence has been transformed into forward and backward sweep, and a function $\overline{FuncName}$ can be created by $\textit{newFunction}$, then in the end of the backward sweep $\textit{FinalizeEE}$ is added. The $\textit{isEstErrFunction}$ is true when the user annotated a function with $\textit{estimate_error}$. Then the next assignment rules extend the known set of reverse mode AD rules.
% Assignment Stmt - Live, Dependent
\[
 \inference[\texttt{\rsema{sema:Assign}:}]
%\inference[]
{
  \textcolor{OliveGreen}{\textit{isLive}}(Stmt) \quad
  \textcolor{OliveGreen}{\textit{isDiff}}(Ref) \\
  Stmt -> Ref \;\textcolor{Blue}{\texttt{=}}\; Expr \\
  \textcolor{Blue}{\texttt{typeof}}(Ref) \xrightarrow {newLocal} \overline{Var} \\
  Expr, \overline{Var} \xrightarrow {Expr} \overline{Stmts} \\
  Ref, Expr, \overline{Var} \xrightarrow {buildAssignError} Est
}
{
  \makecell[l]{
    Ref \;\textcolor{Blue}{\texttt{=}}\; Expr \\
    \quad ->
      \left[
        \makecell[l]{
          \textcolor{OliveGreen}{\texttt{Push}}\textcolor{Blue}{(}Ref\textcolor{Blue}{);}\; Ref \;\textcolor{Blue}{\texttt{=}}\; Expr;\; \overline{Var} \;\textcolor{Blue}{\texttt{=}}\; 0 \\
          \textcolor{OliveGreen}{\texttt{Pop}}\textcolor{Blue}{(}Ref\textcolor{Blue}{);}\;\overline{Stmts} \textcolor{Blue}{;}\;Est
        }
      \right.
  }
}
\]

In rule S2, %TODO: Use ref S\ref{sema:Assign}
the premises are fulfilled when:
\begin{enumerate*}[label=(\roman*)]
    \item $\textit{isLive}$ is true when the statement $\textit{Stmt}$ is useful for the derivative computation;
    \item $\textit{isDiff}$ is true when the memory location is relevant for the derivative computation;
    \item $\textit{newLocal}$ can create a variable of the type of $\textit{Ref}$;
    \item $\textit{Expr}$ can differentiate $\textit{Expr}$; and
    \item \textit{buildAssignError} can generate EE instructions \textit{Est}.
\end{enumerate*}
When all antecedents are fulfilled, then $\textit{AssignError}$ augments the backward sweep. \textit{Ref} represents a memory location where a value can be stored. 

%\textit{Ref} represents a memory location where a value can be stored. The $\textit{isLive}$ is true when the statement $\textit{Stmt}$ is useful for the derivative computation. The predicate $\textit{isDiff}$ is true when the memory location is relevant for the derivative computation. When both predicates are fulfilled and variables of the type of $\textit{Ref}$ can be created by $\textit{newLocal}$ and $\textit{Expr}$ can differentiate $\textit{Expr}$ and \textit{buildAssignError} can generate EE instructions \textit{Est}, then $\textit{AssignError}$ augments the backward sweep.
% Assignment Stmt - Live, Independent
\[
 \inference[\texttt{\rsema{sema:AssignLive}:}]
%\inference[]
{
  \textcolor{OliveGreen}{\textit{isLive}}(Stmt) \quad
  \neg \textcolor{OliveGreen}{\textit{isDiff}}(Ref) \\
  Stmt -> Ref \;\textcolor{Blue}{\texttt{=}}\; Expr \\
  Ref, Expr \xrightarrow {buildAssignError} Est}
{
  Ref \;\textcolor{Blue}{\texttt{=}}\; Expr ->
    \left[
      \makecell[l]{
        \textcolor{OliveGreen}{\texttt{Push}}\textcolor{Blue}{(}Ref\textcolor{Blue}{);}\; Ref \;\textcolor{Blue}{\texttt{=}}\; Expr \\
        \textcolor{OliveGreen}{\texttt{Pop}}\textcolor{Blue}{(}Ref  \textcolor{Blue}{);} \; Est
      }
    \right.
}
\]
\vskip 4pt
% Assignment Stmt - not Live, Dependent
\[
\inference[\texttt{\rsema{sema:AssignNotLive}:}]
%\inference[]
{ 
  \neg \textcolor{OliveGreen}{\textit{isLive}}(Stmt) \quad
  \textcolor{OliveGreen}{\textit{isDiff}}(Ref) \\
  Stmt -> Ref \;\textcolor{Blue}{\texttt{=}}\; Expr \\
  \textcolor{Blue}{\texttt{typeof}}(Ref) \xrightarrow {newLocal} \overline{Var} \\
  Expr, \overline{Var} \xrightarrow {Expr} \overline{Stmts} \\
  Ref, Expr, \overline{Var} \xrightarrow {buildAssignError} Est 
}
{
  Ref \;\textcolor{Blue}{\texttt{=}}\; Expr ->
    \left[
      \makecell[l]{
        \overline{Var} \;\textcolor{Blue}{\texttt{=}}\; 0 \\
        \overline{Stmts} \textcolor{Blue}{;} \; Est
      }
    \right.
}
\]
The next two rules S3 and S4 %TODO: Use ref S\ref{sema:AssignLive} and S\ref{sema:AssignNotLive}
consider the variations of the values of the $\textit{isLive}$ and $\textit{isDiff}$ predicates.

When either $\textit{isDiff}$ or $\textit{isLive}$ is false, then we still create the adjoints and insert the $\textit{AssignError}$ to capture variable's error contribution. Function calls and parameter passing can be expressed by analogy.

%\twocolumn

%% End Semantic description

\subsection{Programming Model \& Framework Design}
\begin{figure}[h]
    \centering
    \includegraphics[width=\columnwidth]{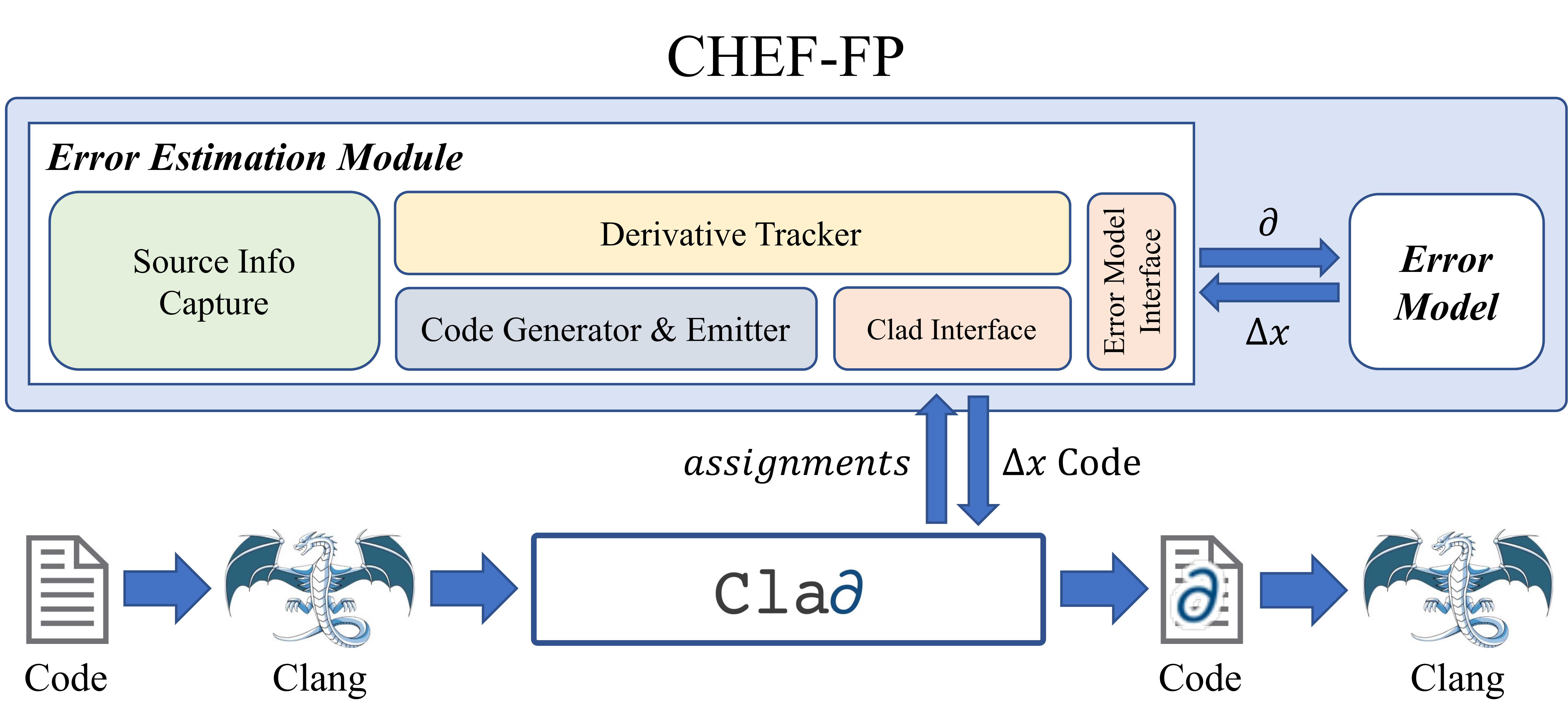}
    \caption{\textbf{FP EE Generation Workflow in \cladtool}. The control from Clang's compilation pipeline is intercepted by Clad, and detection of calls to \mintinline{cpp}{estimate_error()} instantiates the \emph{Error Estimation Module} and sets up callbacks. The \emph{Error Estimation Module} then listens to callbacks from Clad's adjoint mode and is responsible for augmenting the derivative body with the error estimation code defined in the \emph{Error Model} before passing the control back to Clad. The \emph{Error Estimation Module} is also responsible for caching values, tracking derivatives, capturing the required structural and source information, and accumulating errors. }
    \label{fig:CladEEWorkflow}
\end{figure}

Clad's callback system allows the creation of extensions that can augment  generated code. We use this ability to build a lightweight framework to insert FP error estimation code in Clad generated adjoints. CHEF-FP leverages the flexible design of Clad and adds itself as a native extension that synthesizes error estimation code as part of the differentiation process.  Fig.~\ref{fig:CladEEWorkflow} outlines its high-level design. \cladtool's implementation is broadly  divided into an \emph{Error Estimation Module} and an \emph{Error Model}.

\textbf{Error Estimation Module:}
The error estimation module is the major component that makes up CHEF-FP. The registration of calls to \mintinline{cpp}{estimate_error()}, by which the user annotates the functions of interest, instantiates the  \emph{Error Estimation Module}. This sets up error estimation callbacks and passes the control back to Clad. The \emph{Error Estimation Module} listens to callbacks from Clad's adjoint mode and takes control when an interesting (from EE perspective) callback is triggered. The \emph{Error Estimation Module} augments the derivative body by adding the EE expressions received from the \emph{Error Model} and passes control back to Clad until the next callback is invoked. Thus, it is responsible for interfacing with Clad, defining the various methods that facilitate the generation and emission of FP error estimation code, and exchanging information with the \textit{Error Model} to produce specifically augmented code. The error estimation module also caches values, tracks derivatives, captures the required structural and source information, and accumulates errors to provide the final FP error estimate.

\textbf{Error Model:}
The Error Model defines an interface to describe the error expression to generate. It receives derivative expressions as the derivative is being generated and returns respective error expressions. By default, a call to a user-defined function \mintinline{cpp}{getErrorVal} is synthesized, allowing users to write their formulas in plain C/C++ code. An additional interface that can be programmed to generate advanced calls that export more adjoint information to incorporate other advanced custom \emph{Error Model}s also exists. This interface can be completely customized to make the FP error analysis more specific to an application. Listing~\ref{lst:minimalExmp} demonstrates the programming model with a minimal example.

\begin{listing}[h]
\begin{minted}[fontsize=\small, breaklines]{cpp}
float func(float x, float y) {
  float z;
  z = x + y;
  return z;
}
int main() {
  // Call estimate error on target function.
  auto df = clad::estimate_error(func);
  // Declare the inputs, their derivative 
  // outputs and the final error output.
  float x = 1.95e-5, y = 1.37e-7;
  float dx = 0, dy = 0;
  double fp_error = 0;
  // Execute the generated code.
  df.execute(x, y, &dx, &dy, fp_error);
  // fp_error now contains the error of func.
  std::cout << "Error in func: " << fp_error;
}
\end{minted}
\caption{Minimal demonstrator of the usage of CHEF-FP.}
\label{lst:minimalExmp}
\end{listing}

%\begin{figure}[h]
%    \centering
%    % \includegraphics[width=\linewidth]{CladEEWorkflow.pdf}
%    \includegraphics[width=\linewidth]{cladEEWorkflow.png}
%    \caption{FP Error Estimation Generation Workflow in Clad}
%    \label{fig:CladEEWorkflow}
%\end{figure}

\subsection{Implementation}

\cladtool registers an API (\mintinline{cpp}{estimate_error}) in Clad that can be used to calculate the floating point errors in a given function using a default error estimation model. An example of a typical invocation of \cladtool is illustrated in listing~\ref{lst:minimalExmp}.

For more complex analyses, it may be necessary to change the underlying error model to achieve satisfactory estimates. This can simply be achieved by implementing the \mintinline{cpp}{FPErrorEstimationModel} interface with the appropriate error model, compiling it into a shared library, and passing that library to \cladtool. An example of such an implementation is described in listing~\ref{lst:custClass}.
\begin{listing}[h]
\begin{minted}[fontsize=\small, breaklines]{cpp}
struct CustomModel : public FPErrorEstimationModel {
  CustomModel(DerivativeBuilder& builder)
      : FPErrorEstimationModel(builder) {}
  // Returns the error expression to be 
  // calculated for each variable assignment.
  clang::Expr* AssignError(StmtDiff refExpr, const char* name) override;
};
\end{minted}
\caption{A template for declaring custom model classes.}
\label{lst:custClass}
\end{listing}
Here \mintinline{cpp}{AssignError} builds the code expression (expressed as Clang's internal expression types) that can then be emitted into code by the error estimation module. This function exposes three values that can be used to implement the custom model - the name of the variable, the variable itself, and its derivative.

% \mintinline{cpp}{FPErrorEstimationModel} defines another function \mintinline{cpp}{SetError} which is used to initialize error accumulation variables. This function can be either explicitly implemented to initialize variable errors with some specific values or can be left unimplemented, in which case the error accumulation variables initialized with 0.

\begin{listing}[h]
\begin{minted}[fontsize=\small, breaklines]{cpp}
namespace clad { 
  double getErrorVal(double dx, double x, const char* name) {
    return dx * (x - (float)x);
  }
}
clang::Expr* CustomModel::AssignError(StmtDiff refExpr, std::string name) {
  // Build a vector-like container to store
  // the parameters of the function call.
  llvm::SmallVector<clang::Expr*, 3> params{
      refExpr.getExpr_dx(), refExpr.getExpr(), utils::CreateStringLiteral( m_Context, name)};
  // Return a call to getErrorVal.
  return GetFunctionCall("getErrorVal", "clad", params);
}
\end{minted}
\caption{\textbf{An example implementation of \mintinline{cpp}{AssignError} that builds calls to external functions}. The  implementation generates a call to a user-defined function \mintinline{cpp}{getErrorVal}.}
\label{lst:funcAssign}
\end{listing}
There are two broad ways to implement \mintinline{cpp}{AssignError} - building an arithmetic expression or building an external function call. Both methods involve building an assignable expression. The former approach involves directly building the arithmetic expression to be assigned and so it is fairly limited. It also requires that the model be re-built for every modification in \mintinline{cpp}{AssignError}. While we provide an API useful for building simple expressions, it becomes unintuitive to implement complex models based on just a single expression. Hence, we build calls to external functions as a valid error model as long as the function has a compatible return type to the variable being assigned the error. This approach allows users to define their error models as regular C++ functions, allowing for the implementation of more computationally complicated models. 

Listing~\ref{lst:funcAssign} demonstrates how to implement an arbitrary function \mintinline{cpp}{getErrorVal} and use it as a custom error model. In the listing, we build the model used in ADAPT-FP, whose mathematical notation is described as follows:
\begin{equation}
    \label{adapt-model}
\Delta = \sum_{i=1}^{n}{\frac{\delta f}{\delta x_i} * (x_i - (float)x_i)},
\end{equation}
where \(\Delta\) is the accumulated error in the function \(f\) due to the rounding errors for all \(x_i\). Here, the error assigned to each variable is
the difference between its single and double precision values. Note, this model imposes a requirement that all functions it is used on are of double or higher precision. 

It is also possible to modify the signature and composition of the example external error function (i.e., \mintinline{cpp}{getErrorVal}) defined in listing~\ref{lst:funcAssign} by modifying the \mintinline{cpp}{AssignError} function. This allows users to add relevant meta information of the variables exposed by Clang. For example, users can implement \mintinline{cpp}{AssignError} to build calls to an error function that can take in a variable's source location. This information can be used for identifying specific areas of code more prone to errors. 
\section {Experiments}
\label{sec:results}

We compare \cladtool against the current state-of-the-art ADAPT using five different algorithms: \textit{Arc Length}, \textit{Simpsons}, \textit{k-Means clustering}, \textit{HPCCG}, and \textit{Black-Scholes}. For the first 4 benchmarks, we use the error model described in equation \ref{adapt-model}. The time taken is measured using \textit{Google benchmark}
% ~\cite{googlebench} 
and peak memory by \textit{GNU time}. The benchmarks were done on the Princeton \textit{Tiger} cluster~\cite{princetontiger} with a 2.4GHz Intel Xeon Gold 6148 CPU and 188 GB of RAM. 
Similar benchmark results were obtained using the LLNL Quartz~\cite{llnlquartz} system, which is a cluster consisting of Intel Xeon E5-2695 processors with 2.1 GHz cores and 128 GB of memory per node.
% \mynote{Harshitha: can you mention a bit more details on the machines you ran}

\begin{table}[h]
\centering
\def\arraystretch{1.5}
\begin{tabular}{ccccc}
\toprule
\textbf{Benchmark} & \textbf{Threshold} & \textbf{Actual Error} & \textbf{\def\arraystretch{1}\begin{tabular}[c]{@{}c@{}}Estimated \\ Error\end{tabular}} & \textbf{Speedup} \\ \midrule
Arc Length         & 1e-05              & 3.24e-06              & 3.24e-06                                                            & 1.11             \\
Simpsons           & 1e-06              & 7.80e-08              & 1.32e-07                                                             & 2.25             \\
$k$-Means            & 1e-06              & 0.00e+00                    & 0.00e+00                                                                   & -                \\
HPCCG              & 1e-10              & 5.21e-12              & 5.92e-11                                                            & 1.08             \\ \bottomrule
\end{tabular}
\caption{\textbf{Error and performance measurements of the mixed precision versions of the benchmarks}. The table shows a comparison between the actual error and \cladtool's estimated error in the mixed precision versions of the original program. It also shows the execution speedup of the mixed variant. For \emph{$k$-Means}, \cladtool's identified mixed precision configuration for the defined threshold showed no speedup.}
\label{tab:speedups}
\end{table}
\def \figureColumnsWidth {0.32} % 0.48 - 2 columns; 0.32 - 3 columns
\begin{figure*}
\centering
\begin{minipage}[t]{0.7\textwidth}
% \begin{figure}[H]
    \centering
    \includegraphics[width=\linewidth]{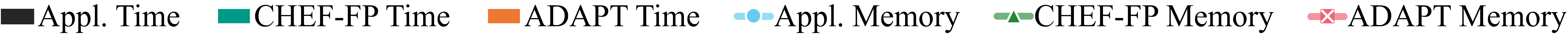}
%    \caption*{Legend for Fig.~\ref{fig:ArclenAvCBench} to Fig.~\ref{fig:BlkSolAvCBench}}
% \end{figure}
\end{minipage}%
\quad
\begin{minipage}[t]{\figureColumnsWidth\textwidth}
%\begin{figure}[h]
    %\centering
    \includegraphics[width=\linewidth]{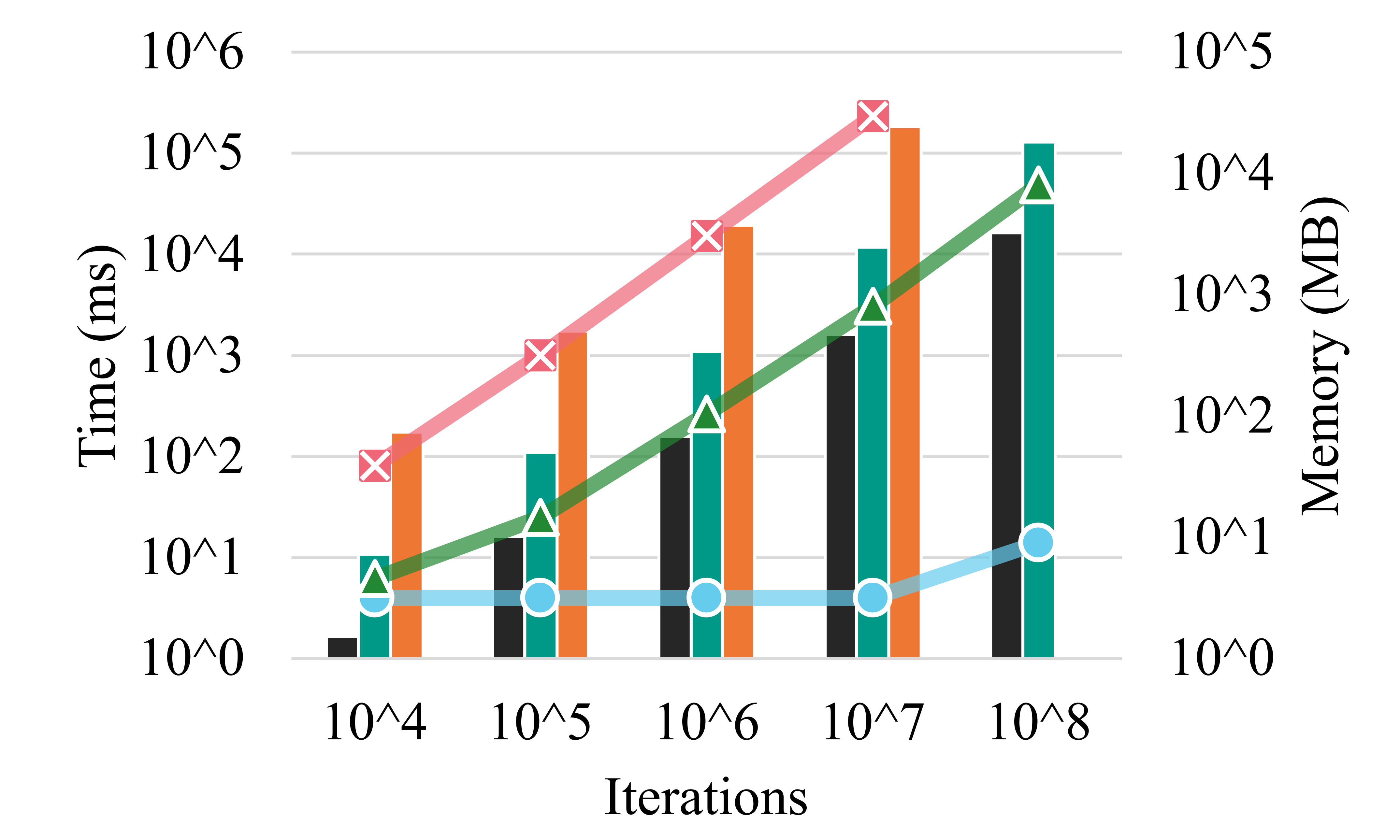}
    \caption{Arc Length}
    \label{fig:ArclenAvCBench}
%\end{figure}
\end{minipage}%
\quad
\begin{minipage}[t]{\figureColumnsWidth\textwidth}
%\begin{figure}[h]
    %\centering
    \includegraphics[width=\linewidth]{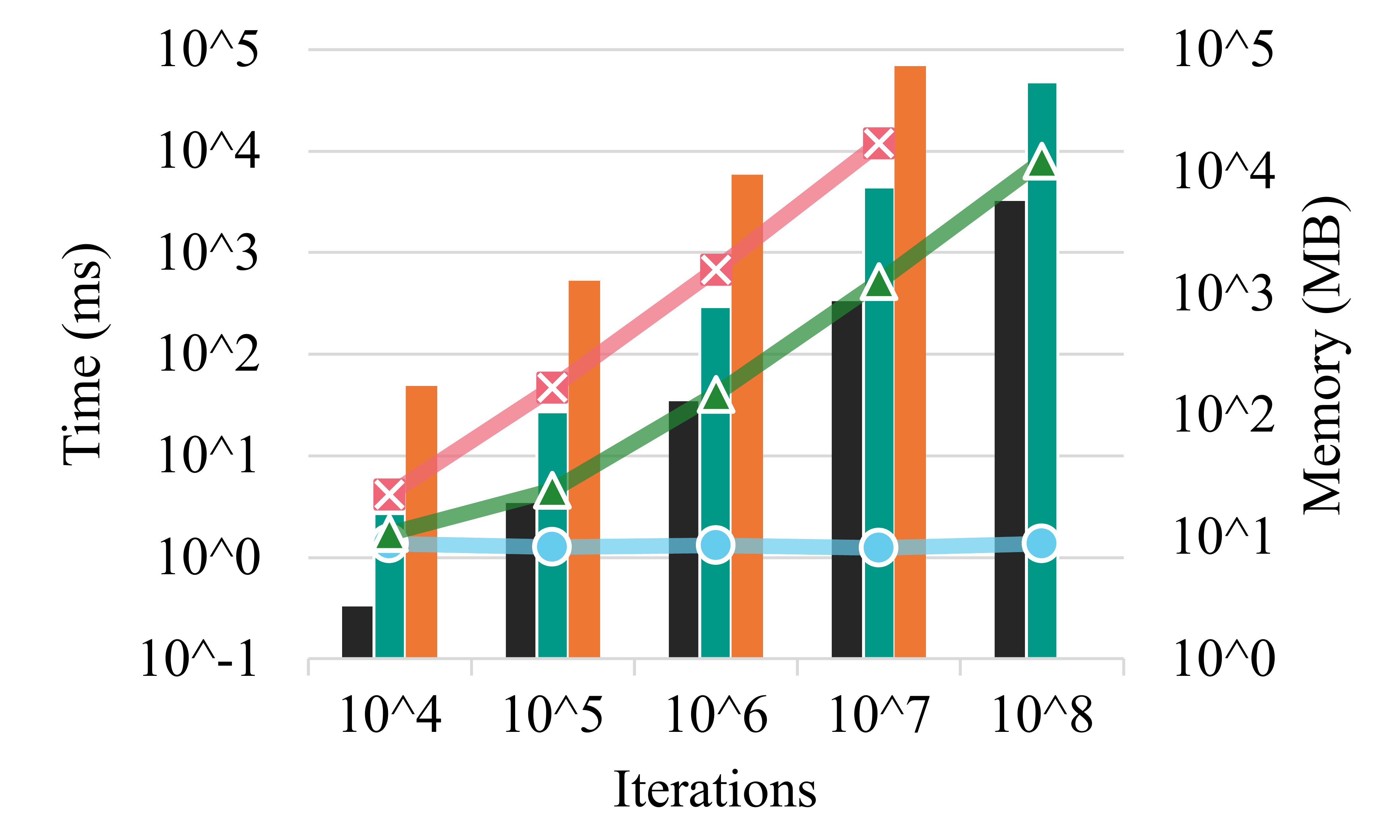}
    \caption{Simpsons}
    \label{fig:SimpAvCBench}
%\end{figure}
\end{minipage}%
\quad
\begin{minipage}[t]{\figureColumnsWidth\textwidth}
%\begin{figure}[h]
    %\centering
    \includegraphics[width=\linewidth]{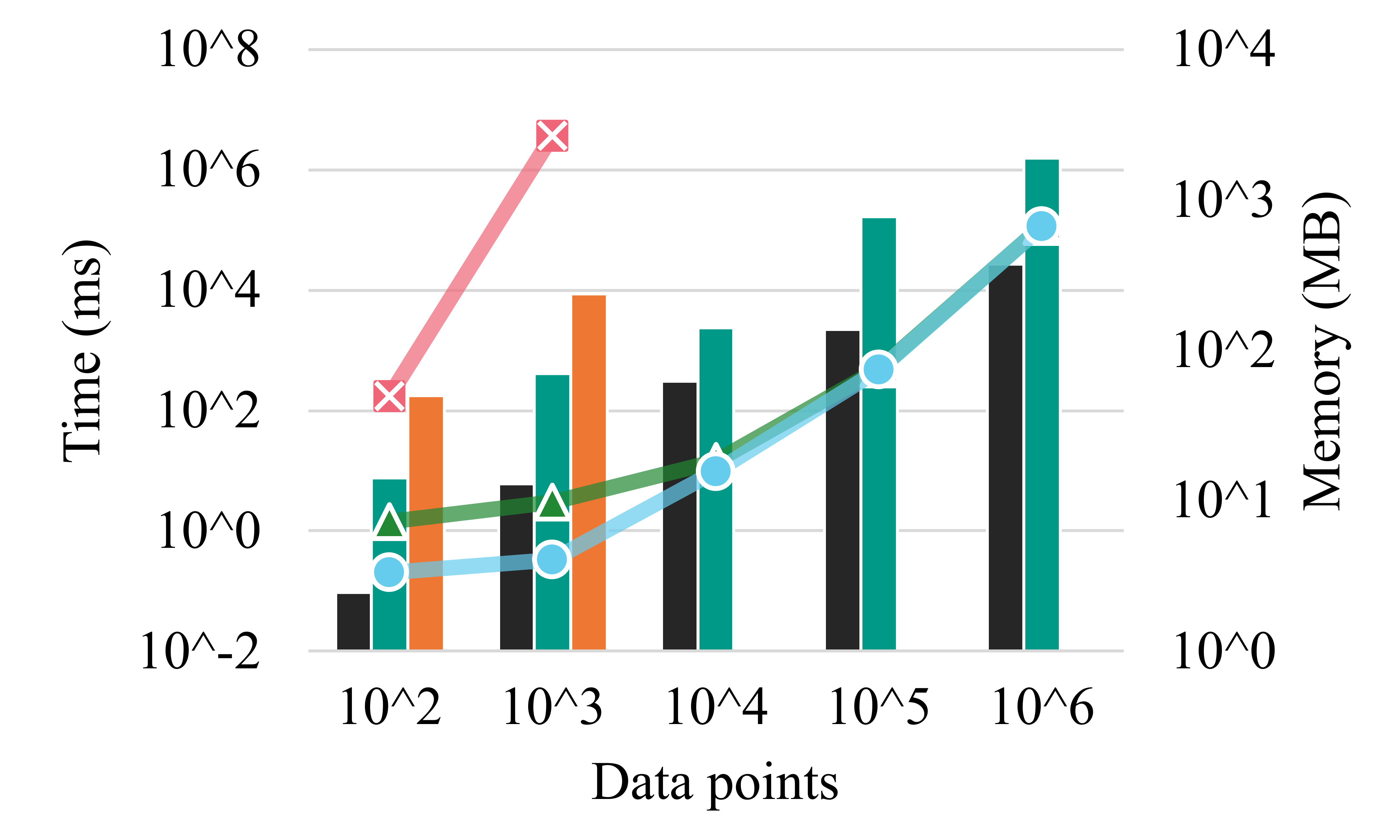}
    \caption{$k$-Means algorithm}
    \label{fig:KMeansAvCBench}
%\end{figure}
\end{minipage}%
\quad
\vskip 12pt
\begin{minipage}[t]{\figureColumnsWidth\textwidth}
%\begin{figure}[h]
    %\centering
    \includegraphics[width=\linewidth]{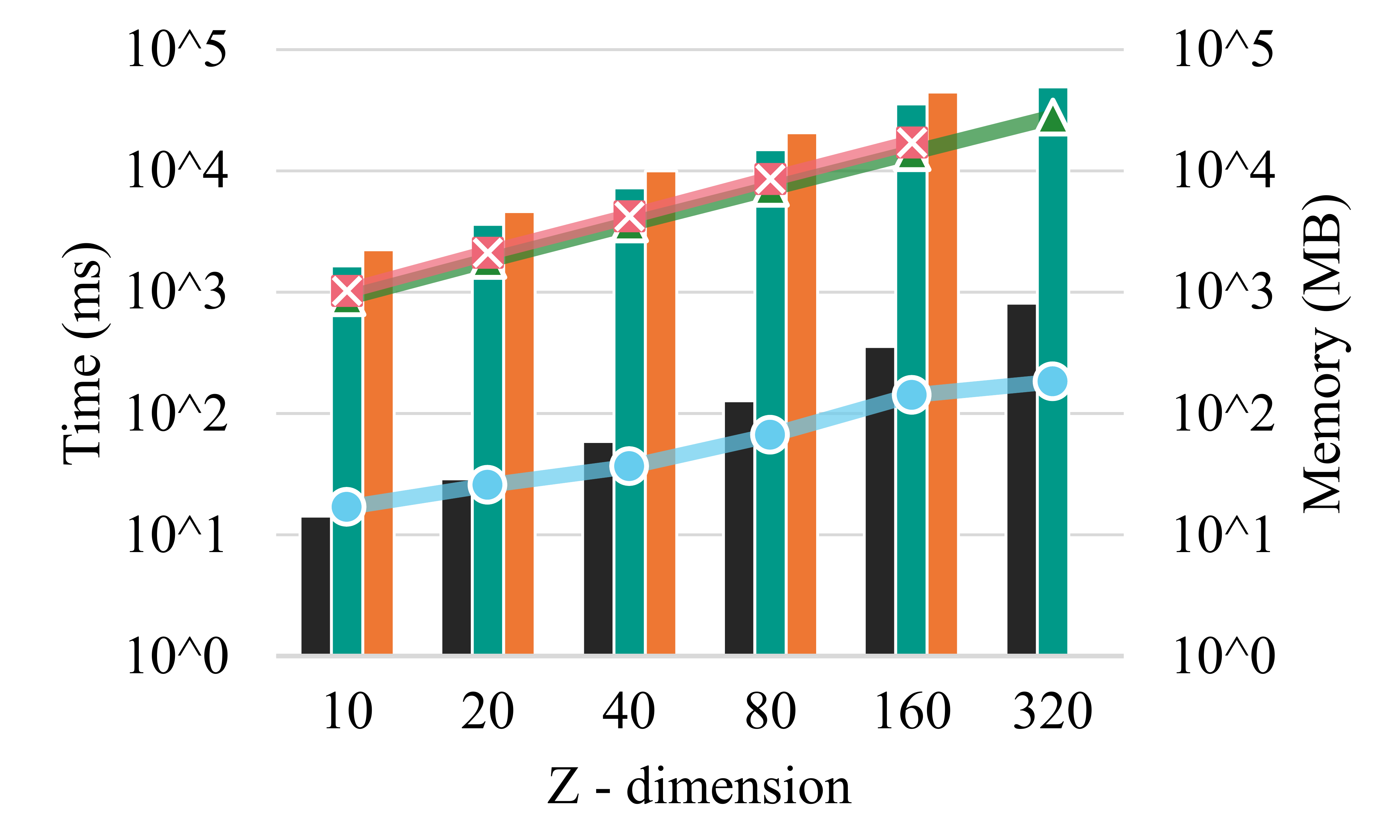}
    \caption{HPCCG}
    \label{fig:HPCCGAvCBench}
%\end{figure}
\end{minipage}%
\quad
\begin{minipage}[t]{\figureColumnsWidth\textwidth}
%\begin{figure}[h]
    %\centering
    \includegraphics[width=\linewidth]{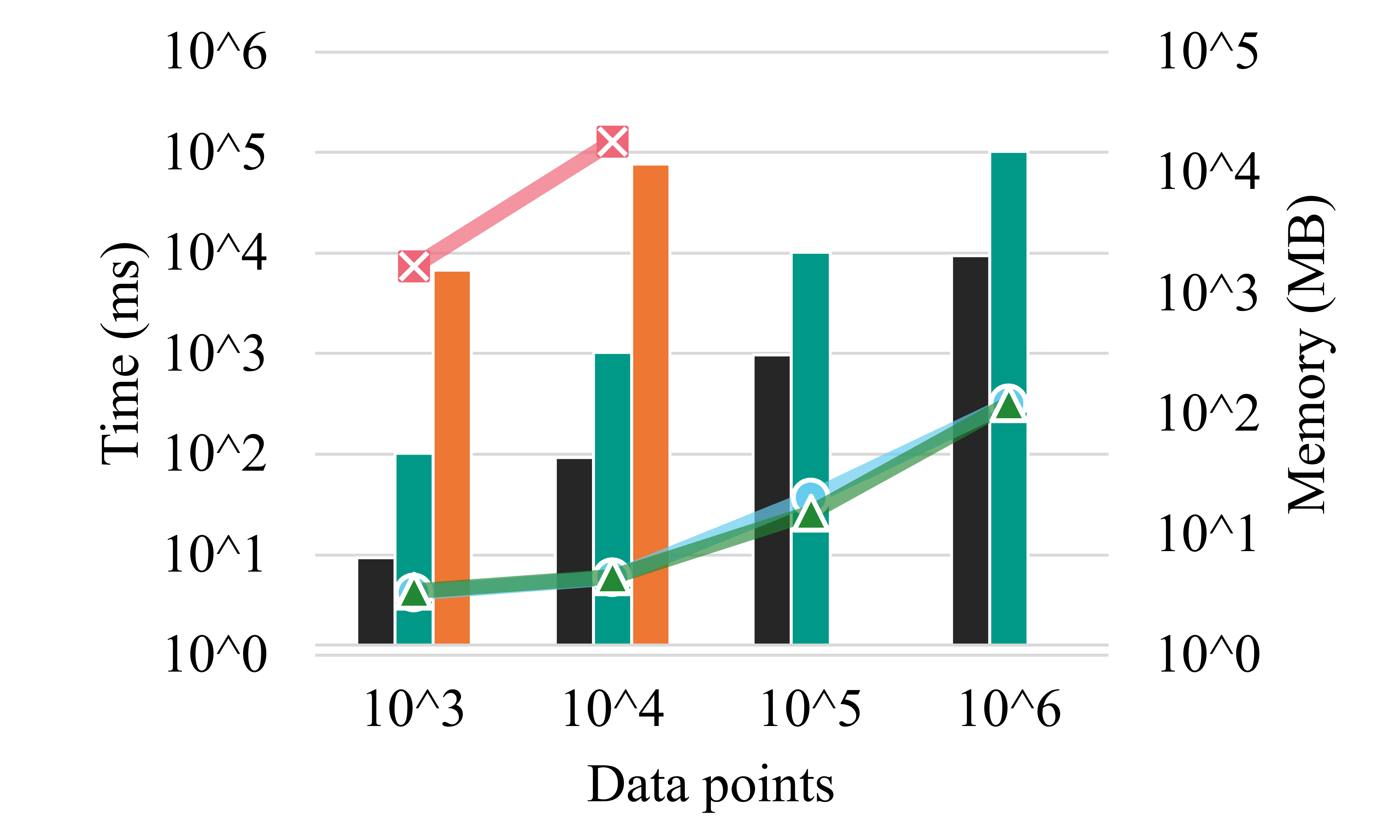}
    \caption{Black-Scholes}
    \label{fig:BlkSolAvCBench}
%\end{figure}
\end{minipage}%
\quad
%The following alignment is a BIG HACK, DO NOT CHANGE PLS
\begin{minipage}[t][0.60cm][b]{\figureColumnsWidth\textwidth}
    \def\arraystretch{1.35}
    \centering
    \setlength\tabcolsep{3pt}
    \resizebox{!}{0.82\height}{%
    \begin{tabular}{ccc}
    \toprule
    \textbf{Benchmark} & \textbf{Time} & \textbf{Memory} \\ \midrule
    Arc length           & 1.61x             & 1.95x                     \\
    Simpsons             & 2.17x             & 1.44x                     \\
    $k$-Means            & 2.02x             & 4.44x                     \\
    HPCCG                & 1.03x             & 1.02x                     \\
    Black-Scholes    & 1.76x           & 6.32x                           \\ \bottomrule
    \end{tabular}
    }
    \vspace{4pt}
    \captionof{table}{Performance Improvements}
    \label{tab:chefperfsum}

\end{minipage}%

\quad
\caption*{\textbf{Figures  \ref{fig:ArclenAvCBench}
to
\ref{fig:BlkSolAvCBench}
show the results of benchmarking \cladtool, ADAPT, and the original function.} The labels below the function show the algorithm being benchmarked. The lines show peak memory usage and the bars represent the time taken during the FP error analysis of the given algorithm. \textbf{Table~\ref{tab:chefperfsum} summarises CHEF-FP's performance improvements over ADAPT.} The improvements are given as 'times improved' over the FP error analysis and represent the average improvement across all the data points. Our benchmarks show that CHEF-FP outperforms ADAPT, which is the current state-of-the-art AD-based FP error estimation tool while producing mixed precision analysis results that agree with ADAPT's analysis.}
\vspace{-0.5cm}
\end{figure*}
For all of the following benchmarks, we present reductions in the analysis time and memory while still producing results that are at par with tools such as ADAPT. We also utilize this section to comment on how the analysis for certain examples can be extended and how \cladtool can be used to perform more sophisticated approximate analysis. Lastly, we present a summary of the mixed precision analysis results in table~\ref{tab:speedups} to demonstrate how \cladtool can produce correct floating-point analysis results while requiring little to no user intervention.

% We also report speedups from implementing \cladtool's mixed precision recommendations and outline the estimated errors of the mixed configurations when compared to the original high-precision version of the same application. These results are given in table .

% \begin{table}[htbp]
% \caption{Table Type Styles}
% \begin{center}
% \begin{tabular}{|c|c|c|c|}
% \hline
% \textbf{Table}&\multicolumn{3}{|c|}{\textbf{Table Column Head}} \\
% \cline{2-4} 
% \textbf{Head} & \textbf{\textit{Table column subhead}}& \textbf{\textit{Subhead}}& \textbf{\textit{Subhead}} \\
% \hline
% copy& More table copy$^{\mathrm{a}}$& &  \\
% \hline
% \multicolumn{4}{l}{$^{\mathrm{a}}$Sample of a Table footnote.}
% \end{tabular}
% \label{tab1}
% \end{center}
% \end{table}

\subsubsection{Arc Length}

The arc length function approximates a curve's length ($L$) by sampling various points on the curve and summing up the straight line distance between two consecutive points on the curve.
\[ L = \lim_{n\to\infty} \sum_{i=1}^{n} \sqrt{{\Delta x}^2 + {\Delta y_i}^2}
\]

We vary the number of iterations, $n$, that the arclength algorithm is run for, to benchmark \cladtool against ADAPT. Fig.~\ref{fig:ArclenAvCBench} compares the time taken and memory used to analyze the algorithm for mixed precision analysis. The absence of a data point for ADAPT at $10^8$ iterations is due to it running out of memory. \cladtool's memory footprint is lower because the error calculation is inlined in the gradient function. 

\subsubsection{Simpsons}

Simpsons is an iterative algorithm to approximate the integral of a function in the given interval by summing the integral over multiple small intervals:
\[ \int_{a}^{b} f(x)\, dx \approx \frac{h}{3} \left[f(a) + f(b) + 4 \sum_{i = 1, 3, 5}^{2n - 1} f_i + 2 \sum_{i = 2, 4, 6}^{2n - 2} f_i
\right]
\]
Here, \(f_i = f(a+ih)\), \(h = (a+b)/2n\), and \(n\) is the number of iterations. Fig.~\ref{fig:SimpAvCBench} compares ADAPT and \cladtool on the basis of time and memory usage. Similar to arc length we vary the number of iterations to benchmark the two tools. \cladtool's recommended mixed precision configuration gives a speedup of 2.25 times when compared against the same program in higher precision. It is also able to predict the actual error in the mixed-precision version of the application, as seen in table~\ref{tab:speedups}.

\subsubsection{$k$-Means Clustering}

Part of the Rodinia benchmark suite~\cite{5306797}, the $k$-Means clustering algorithm is used for grouping multiple 
data points into $k$ clusters. We instrument the Euclidean distance function as it is the major computational hotspot of the application. Similar to previous benchmarks, we compare the performance of \cladtool against ADAPT in fig.~\ref{fig:KMeansAvCBench}. 

The Euclidean distance function has three major variables: \textit{attributes}, \textit{clusters}, and \textit{sum}. It can be represented as follows:
\[\textit{sum} = \sqrt{\sum_{i=0}^{n}{(\textit{attributes}_i - \textit{clusters}_i)^2}}\]
\begin{table}[h]
\vspace{2.5mm}
\def\arraystretch{1.5}
\centering
\begin{tabular}{ccc}
\toprule
\textbf{Variable(s) in Lower Precision} & \textbf{Actual Error} & \textbf{Estimated Error} \\ \midrule
\textit{attributes} & 00e+00   & 00e+00   \\
\textit{clusters}  & 3.67e-04 & 9.35e-04 \\
\textit{sum}        & 8.33e-04 & 7.08e-03 \\
all 3               & 2.40e-03 & 8.01e-03 \\ \bottomrule
\end{tabular}
\caption{\textbf{k-Means -- Error measurements of various mixed precision configurations}. We demote the 3 variables (\emph{attributes}, \emph{clusters} and \emph{sum}) to lower precision one by one and compare the resulting errors with \cladtool's estimate.}
\label{tab:kmeans}
\vspace{-3.8mm}
\end{table}
The error estimated by Clad for \textit{attributes} is 0 because the input data of the benchmark is represented with four digits after the decimal. The errors estimated for \textit{clusters} and \textit{sum} are higher than the threshold set for $k$-Means in table~\ref{tab:speedups}. Hence, \cladtool recommends only converting \textit{attributes} to lower precision. To test out \cladtool's error estimates, we went a step ahead and converted each of the variables to lower precision individually and found the actual error introduced by them. These configurations were executed on $10^6$ datapoints, and the findings are shown in table~\ref{tab:kmeans}.

\subsubsection{HPCCG}

Part of the Mantevo benchmark suite, 
HPCGG is a simple conjugate gradient benchmark code for a 3D chimney domain converted to be single-threaded.
The analysis is done while scaling the inputs from the base dimension of $20\times30\times10$ to the recommended size of $20\times30\times160$ and then further to $20\times30\times320$. The results in fig.~\ref{fig:HPCCGAvCBench} show that ADAPT runs out of memory for the $20\times30\times320$ dimension.

\begin{figure}[h]
    \centering
    \includegraphics[width=0.8\columnwidth]{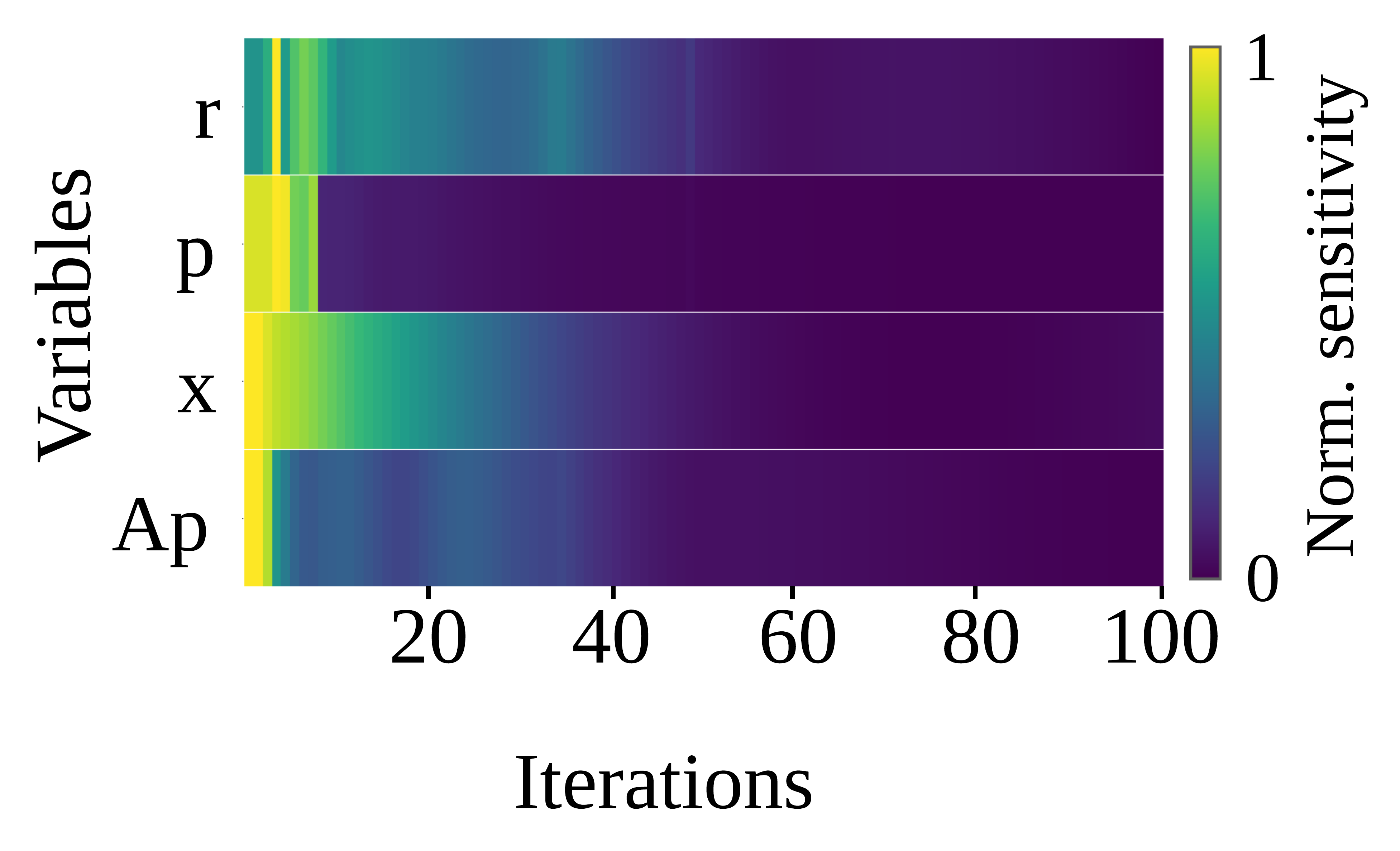}
    \caption{\textbf{HPCCG Variable Heatmap.} This illustrates the normalized sensitivity of the variables \textit{r}, \textit{p}, \textit{x}, and \textit{Ap} for every iteration. This sensitivity profile can be used to determine more fine-tuned optimizations.}
    \label{fig:HPCCGHeatmap}
\vspace{-3.8mm}
\end{figure}

\cladtool is also used to analyze HPCCG for a possible loop perforation based optimization. We already generated errors for intermediate variables, so we only needed to slightly tweak it to dump the sensitivity of the variables over each iteration. We analyzed the change in sensitivity of the various variables over the total run of the application, and we found that the sensitivity for all variables drops below our set threshold after almost $60$ iterations. The normalized sensitivity of the variables is shown in the form of a heat map in fig.~\ref{fig:HPCCGHeatmap}. Based on these findings, we split the main loop of HPCCG into two chunks - the first chunk runs the full loop for the first 60 iterations in high precision, and the second chunk runs for the remaining iterations in lower precision. This configuration gives us a speedup of $8\%$ as shown in table~\ref{tab:speedups}.

\subsubsection{Black-Scholes}

Part of the Parsec-benchmark suite \cite{bienia11benchmarking}, the Black-Scholes equation is a differential equation that describes the change in the value of an option as the price of the  underlying asset changes. 
%It assumes that the asset in consideration follows the Black-Scholes model.
To compare against ADAPT, we analyze the function that calculates the option price depending on the factors that influence the stock. ADAPT is not able to scale beyond $10^4$ data points in our benchmarks. The results are shown in fig.~\ref{fig:BlkSolAvCBench}. For this benchmark, both \cladtool and ADAPT showed that most of the intermediate variables are very sensitive to errors and  could not identify an effective mixed precision configuration.

In addition, we leveraged the customizability of \cladtool to analyze the effect of using Paul Minero's FastApprox library~\cite{fastapprox} instead of the standard C math library to gain performance improvements. The FastApprox library provides approximate versions of various math functions that trade off accuracy for performance. \cladtool can determine the error introduced by replacing the standard versions of the math functions with the approximate ones in any program. We identified three math functions in the Black-Scholes application that had approximate versions in the FastApprox library. The inputs to these functions were identified, and a map was formed with their names mapping to the functions they were input to. The map is used in a custom model to match each variable with the correct function and its approximation. These values are then used by the custom model to accurately estimate the errors due to approximation in the application. The algorithm for the custom model is illustrated in algorithm \ref{lst:algoBlkSch}. The estimated error, actual error, and speedup from FastApprox in the Black-Scholes application are shown in table~\ref{tab:BSFastApprox}.

% \begin{algorithm} 
%  \caption{Error model for estimating approximation-based errors.}
% 	\begin{algorithmic}[1]
% 	\Require $x$ is variable of interest, $dx$ is its partial derivative wrt. the function, and $name$
% 	\ForAll{$n$ in \{``sqrt", ``exp",...\}}
%       \If{$name = ``clad\_"+n$}
%           \State $delta \gets \Call{eval}{n,x} - \Call{eval}{Approx::n,x}$
%       \EndIf
%     \EndFor
%     \State $xApproxError \gets | dx * delta |$
%     \State $\Call{registerError}{name, xApproxError}$
%     \State \Return $xApproxError$
% 	\end{algorithmic}
%  \label{lst:algoBlkSch}
% \end{algorithm} Garima: Just commented this out because this is not how we do it actually, tried to put another version into perspective below. If it is not up to par, feel free to remove it.

\begin{algorithm}[h] 
 \caption{Estimating approximation-based errors.}
	\begin{algorithmic}[1]
	\Require input variable as $x$ and its name as $name$, the partial derivative of $x$ wrt. the function as $dx$, and a map of variables of interest as $S:name ->$ function name
	\State $\Delta \gets 0$
	\If{$name$ is contained in $S$}
	\State $fName \gets$ \Call{S.getValue}{$name$}
	\State $\Delta \gets$ \Call{eval}{$fName$, $x$}$-$\Call{evalApprox}{$fName$, $x$}
	\EndIf
    \State $xApproxError \gets | dx * \Delta |$
    \State \Call{registerError}{$name$, $xApproxError$}
    \State \Return $xApproxError$
	\end{algorithmic}
 \label{lst:algoBlkSch}
\end{algorithm}

\begin{table*}[]
\centering
\def\arraystretch{1.5}
\begin{tabular}{ccccccccc}
\toprule
\multirow{2}{*}{\textbf{App Configuration}} &
  \multicolumn{3}{c}{\textbf{Actual Errors}} &
   &
  \multicolumn{3}{c}{\textbf{Estimated Errors}} &
  \multirow{2}{*}{\textbf{Speedup}} \\ \cline{2-4} \cline{6-8}
                        & \textit{avg.} & \textit{max.} & \textit{acc.} &  & \textit{avg.} & \textit{max.} & \textit{acc.} &      \\ \midrule
FastApprox w/o Fast exp & 1.16e-04      & 9.62e-05      & 1.16e+01      &  & 9.16e-04      & 6.25e-04      & 9.60e+00       & 1.14 \\
FastApprox w/ Fast exp  & 5.8e-04       & 6.9e-04       & 5.88e+01      &  & 3.5e-03       & 5.0e-03       & 1.07e+02       & 1.65 \\ \bottomrule
\end{tabular}
\caption{\textbf{Black-Scholes - Error and performance analysis of the various FastApprox based configurations}. This table shows the actual error and \cladtool's estimated error of the approximate version of the application. More specifically, it outlines the average, maximum, and accumulated error over 1000 data points. For the first row, the original program uses approximate versions of the \emph{log} and \emph{sqrt} functions. For the second row, the application uses the approximate version of \emph{exp}.}
\vspace{-0.5cm}
\label{tab:BSFastApprox}
\end{table*}
\section{Discussion}

\cladtool shows consistently better performance than ADAPT while suggesting similar hints for mixed precision tuning and yielding similar performance benefits. \cladtool provides an efficient, straightforward, and flexible way to analyze FP errors in complex C++ applications. We demonstrate the use of \cladtool for mixed precision tuning and sensitivity analysis on 4 benchmarks -- \emph{Arc Length}, \emph{Simpsons}, \emph{$k$-Means} and \emph{HPCCG}, and its use for approximation analysis on 1 benchmark -- \emph{Black-Scholes}. In this section, we discuss in more detail the experimental results and current limitations.
%In this section, we present a summary of results from Clad's Floating Point Error Estimation Framework, we also discuss the limitations of our work and indicate possible future directions of research.  

%DL saves space
\subsection{Summary of Results From \cladtool}

\cladtool identified stable mixed precision configurations for the \emph{Arc Length} and \emph{Simpsons} benchmarks for a given threshold; on conversion to these mixed precision configurations, both applications saw a noticeable speedup, summarised in table~\ref{tab:speedups}. For the \emph{$k$-Means} benchmark, CHEF-FP could not identify a mixed precision configuration (with a threshold of \(10^{-6}\)) that resulted in a speedup. However, CHEF-FP provided good upper-bound estimates of the error in the application for different mixed precision configurations (table~\ref{tab:kmeans}).

For the \emph{HPCCG} benchmark, the tool was used to analyze the sensitivity of variables across the main loop. This allowed for us to discover a mixed precision configuration that involved splitting the main loop into performing the first 60 iterations in higher precision and the rest in lower precision. Lastly, we also demonstrated the flexibility of CHEF-FP by leveraging its custom model support to perform an approximation error analysis on the \emph{Black-Scholes} benchmark. \cladtool was able to accurately quantify the errors related to approximation in a set of specific functions. Additionally, \cladtool generated two different approximation configurations, and a report on the estimated errors from the same is shown in table~\ref{tab:BSFastApprox}.

\subsection{Current Limitations}
%\cladtool is an easy-to-use and effective tool in analysing FP errors, it still suffers from some limitations:

\textbf{Quantifying overhead of type-casts}: It is possible for a mixed-precision configuration to show worse performance than the high-precision version. Usually, this is due to the overhead of the implicit type-casts that have been introduced in the code by lowering the precision of some variables. One way to combat these overheads is to keep track of them via counters; for example, a trivial implicit casts counter can be implemented using Clang's AST Matchers. %~\cite{clangASTMatch}

\textbf{Variety of analysis datasets}: The results of the mixed-precision configurations formed here are input dependent. 
%While \cladtool uses a variety of inputs to verify the correctness of the mixed-precision configuration, it does not cover the whole domain of the application being analysed and as such, its analysis is fairly dependent on the data used. 
To form a general mixed-precision configuration, it is important to analyze the application over a representative set of inputs.
%dataset that captures its full range. 

% \textbf{Complexity of writing advanced custom models that return expressions}: While Clad does provide a simple interface for function-call based error models, some prerequisite knowledge of the Clang API is required to to build models that output arithmetic expressions. However, for a lot of the complex Clang APIs, Clad provides simpler, short-hand functions with thorough documentation and demos. \mynote{emph on what is there in clad already.}

% \textbf{Functions in external libraries}: It is not possible for Clad to access the source code of functions in external shared libraries, as such Clad cannot estimate the FP error in those functions. The only way to tackle this is to either expose the source of the selected functions to Clad, or incorporate the expected FP errors of these functions in the custom models.

\textbf{Source rewriting for mixed precision configurations}: 
\cladtool{} provides sensitivity profiles, and error estimates to guide the process of mixed-precision re-implementation. Currently, we manually rewrite the source code to implement the mixed precision configurations suggested by \cladtool{}. We can use source transformation tools, such as Typeforge~\cite{pinnow2019typeforge}, to automate the generation of mixed-precision code.
%Currently, the user has to manually rewrite the source code to implement the mixed precision configurations suggested by \cladtool. This can be especially tedious for complex programs. 
As future work, this process can be automated by combining the decision-making and code generation by using the error information at runtime to just-in-time optimize areas of code with lower sensitivity to run in lower precision.

\textbf{Compiler optimizations:} Certain floating-point unsafe optimizations (such as \textit{{\emph --ffast-math}} or \textit{{\emph -fp-model fast}}) can cause CHEF-FP's predicted errors to be different than the actual errors. Currently, CHEF-FP cannot distinguish errors introduced by optimizations from the expected FP errors. Since these optimizations take place post the derivative generation, certain substitution optimizations may cause even the derivative to be incorrect, leading to the underlying error propagation to also be incorrect. Users of  CHEF-FP should be careful as these optimizations can cause incorrect analysis results. Another way the analysis results can be affected is by changing the intermediate rounding mode for mixed precision expressions (through flags such as \textit{{\emph -fp-model}}). Changing this may cause the mixed precision version of the program to suffer performance degradation. We recommend using the \textit{source} mode for rounding of intermediate calculations for consistent results.
\section{Related work}
Many floating point error analysis tools have been proposed in literature, including both static and dynamic techniques. 
Dynamic approaches require running the program to gather necessary information to perform analysis. 
Brown et al.~\cite{brown2007profiling} designed FloatWatch, built on Valgrind, to determine if floating-point operations can be optimized by using a lower precision representation or fixed-point arithmetic. This is done by tracking the maximum difference between single and double precision computations by performing 
%single and double precision executions 
both
simultaneously. Benz et al.~\cite{benz2012dynamic} presented an approach where every floating-point computation is executed side by side in higher precision to assist the programmer in locating floating-point accuracy problems. 
Lam et al.~\cite{lam2013dynamic} proposed a dynamic approach using a binary analysis tool, DynInst, to detect floating-point cancellations. An et al.~\cite{an2008fpinst} developed a dynamic binary analysis, FPInst,
based on DynInst to compute errors by applying simple error accumulation formulas and tracking the error throughout a program.
% Several static analysis tools have been proposed, including FPTaylor~\cite{solovyev2018rigorous}, SATIRE~\cite{das2020scalable}, Gappa~\cite{daumas2010certification}, Precisa~\cite{feliu2018abstract}, Fluctuat. 
Static analysis tools, such as FPTaylor~\cite{solovyev2018rigorous}, SATIRE~\cite{das2020scalable}, and Precisa~\cite{feliu2018abstract}, use a global optimizer to estimate the upper bound of rounding errors. Gappa~\cite{daumas2010certification} automates error evaluation and propagation using interval arithmetic. Static analysis approaches %have the advantage that they
provide a rigorous error analysis, but have been applied only to small benchmarks.
SEESAW~\cite{das2021robustness} employs symbolic adjoint mode AD to give tighter error bounds for intervals of input values and has been shown to work on practical HPC benchmarks. None of these methods targeted mixed-precision.

Several efforts have evaluated whether a program can take advantage of mixed-precision. Most of the techniques used search-based optimization to select suitable mixed-precision versions of the program that satisfies a user-provided error threshold. Lam et. al~\cite{lam2013automatically, feliu2018abstract} proposed CRAFT which uses a search algorithm to automate the identification of code regions that can use lower precision. Gonzales et al.~\cite{rubio2013precimonious} used delta debugging to narrow the search space for mixed-precision configuration. It was extended to consider groups of variables to further reduce the search space~\cite{guo2018exploiting}. Laguna et. al~\cite{laguna2019gpumixer} proposed GPUMixer, a tool used to tune FP precision on GPU programs with a focus on performance improvements. It uses shadow computations analysis to compute the error introduced by mixed-precision for each kernel and uses a search-based technique to identify the best mixed-precision configuration. These methods  work by identifying a set of variables that can be in single precision while leaving the rest of the variables in double precision. Search-based techniques have the drawback that they require several runs of the application, and exploring the space is extremely time-consuming.

There have been several efforts directed towards analyzing applications for
introducing mixed precision as well as estimating the error
due to reduced precision representation using AD. AD has been
used for estimating the rounding error in numerical algorithms since the early 90s~\cite{iri1991history}, where
partial derivatives given by AD were used to 
obtain a first-order approximation of the global rounding error due to elementary 
rounding errors. Interval analysis with the mean value theorem was used to provide tighter
upper-bound for rounding error estimation. Later Langlois~\cite{langlois2000revised} proposed the CENA method, where a correction term
was introduced to the first-order effect of rounding errors on the output of the numerical algorithms
to improve the accuracy of estimation. Subsequently, ADAPT~\cite{menon2018adapt} used AD to estimate errors, enabling mixed-precision tuning by identifying regions where lower precision can be applied while staying within an error threshold.
While ADAPT provided guidance for mixed-precision implementation, it involved manual annotations
and code transformations. To address this issue,  FloatSmith~\cite{lam2019floatsmith} was introduced, it
integrated ADAPT with Codipack (AD tool)~\cite{SaAlGauTOMS2019} and Typeforge (based on Rose \cite{quinlan2011rose} compiler)
to automate the process of analyzing numerical codes. However, this long toolchain made it slow and cumbersome to work with. 
\section{Conclusion}
In this paper, we have defined formalism to augment automatic differentiation to perform floating-point error analysis, and demonstrated an efficient tool using compiler-based source transformation that
%We have laid out a theoretical approach to floating point error estimation and demonstrated it through an architectural enhancement of a compiler-based approach to automatic differentiation. 
does not overwhelm the already overly complex reverse accumulation AD mode. We present \cladtool, a flexible, scalable, and easy-to-use source-code transformation AD-based tool for the analysis of approximation errors in HPC applications. \cladtool works on the source level to inject error estimation code into generated adjoints, allowing analysis to be sped up via compiler optimizations. This setup allows \cladtool to operate on higher memory loads when compared to other FP error estimation tools. It provides considerable flexibility on what estimation code is generated by using custom error models, facilitating the exploration  of other areas of approximation error analysis. 
%with ease.

We demonstrated that \cladtool performs the same analysis as ADAPT-FP in a time and memory-efficient manner using five benchmarks. At analysis time, \cladtool obtained a maximum speedup of 2.17x over ADAPT-FP for the Simpsons benchmark and a memory reduction of 6.32x for the Black-Scholes benchmark. We showed how \cladtool could be used to perform sensitivity analysis and further provided recommendations on how to perform mixed-precision tuning on various applications. We illustrated how it could be used to accurately evaluate different precision configurations. Finally, we explored estimating approximation-based errors and evaluated the resulting approximate configurations on the \emph{Black-Scholes} benchmark showing a speedup of up to 65\%.

The open-source artifact for this work is available at the DOI: \href{https://doi.org/10.5281/zenodo.7660443}{10.5281/zenodo.7660443}.
%It is a double blind submission
\section{Acknowledgements}
The authors would like to thank the anonymous reviewers for their valuable comments and suggestions.
This project is supported by National Science Foundation under Grant OAC-1931408 and under Cooperative Agreement OAC-1836650. This material is based upon work supported by the U.S. Department of Energy, Office of Science, Office of Advanced Scientific Computing Research and  by the LLNL-LDRD Program under Project No. 20-ERD-043. This work was performed under the auspices of the U.S. Department of Energy by Lawrence Livermore National Laboratory under Contract DE-AC52-07NA27344 (LLNL-CONF-841475).

\bibliographystyle{IEEEtran}
%\bibliography{IEEEabrv,mybibfile}
\bibliography{paper}

\begin{thebibliography}{10}
\providecommand{\url}[1]{#1}
\csname url@samestyle\endcsname
\providecommand{\newblock}{\relax}
\providecommand{\bibinfo}[2]{#2}
\providecommand{\BIBentrySTDinterwordspacing}{\spaceskip=0pt\relax}
\providecommand{\BIBentryALTinterwordstretchfactor}{4}
\providecommand{\BIBentryALTinterwordspacing}{\spaceskip=\fontdimen2\font plus
\BIBentryALTinterwordstretchfactor\fontdimen3\font minus
  \fontdimen4\font\relax}
\providecommand{\BIBforeignlanguage}[2]{{%
\expandafter\ifx\csname l@#1\endcsname\relax
\typeout{** WARNING: IEEEtran.bst: No hyphenation pattern has been}%
\typeout{** loaded for the language `#1'. Using the pattern for}%
\typeout{** the default language instead.}%
\else
\language=\csname l@#1\endcsname
\fi
#2}}
\providecommand{\BIBdecl}{\relax}
\BIBdecl

\bibitem{rubio2013precimonious}
C.~Rubio-Gonz{\'a}lez, C.~Nguyen, H.~D. Nguyen, J.~Demmel, W.~Kahan, K.~Sen,
  D.~H. Bailey, C.~Iancu, and D.~Hough, ``Precimonious: Tuning assistant for
  floating-point precision,'' in \emph{SC'13: Proceedings of the International
  Conference on High Performance Computing, Networking, Storage and
  Analysis}.\hskip 1em plus 0.5em minus 0.4em\relax IEEE, 2013, pp. 1--12.

\bibitem{lam2018fine}
M.~O. Lam and J.~K. Hollingsworth, ``Fine-grained floating-point precision
  analysis,'' \emph{The International Journal of High Performance Computing
  Applications}, vol.~32, no.~2, pp. 231--245, 2018.

\bibitem{das2021robustness}
A.~Das, T.~Tirpankar, G.~Gopalakrishnan, and S.~Krishnamoorthy, ``Robustness
  analysis of loop-free floating-point programs via symbolic automatic
  differentiation,'' in \emph{2021 IEEE International Conference on Cluster
  Computing (CLUSTER)}.\hskip 1em plus 0.5em minus 0.4em\relax IEEE, 2021, pp.
  481--491.

\bibitem{solovyev2018rigorous}
A.~Solovyev, M.~S. Baranowski, I.~Briggs, C.~Jacobsen, Z.~Rakamari{\'c}, and
  G.~Gopalakrishnan, ``Rigorous estimation of floating-point round-off errors
  with symbolic taylor expansions,'' \emph{ACM Trans. on Programming Languages
  and Systems (TOPLAS)}, vol.~41, no.~1, pp. 1--39, 2018.

\bibitem{menon2018adapt}
H.~Menon, M.~O. Lam, D.~Osei-Kuffuor, M.~Schordan, S.~Lloyd, K.~Mohror, and
  J.~Hittinger, ``Adapt: Algorithmic differentiation applied to floating-point
  precision tuning,'' in \emph{SC18: International Conference for High
  Performance Computing, Networking, Storage and Analysis}.\hskip 1em plus
  0.5em minus 0.4em\relax IEEE, 2018, pp. 614--626.

\bibitem{vassiliadis2016towards}
V.~Vassiliadis, J.~Riehme, J.~Deussen, K.~Parasyris, C.~D. Antonopoulos,
  N.~Bellas, S.~Lalis, and U.~Naumann, ``Towards automatic significance
  analysis for approximate computing,'' in \emph{2016 IEEE/ACM Int. Symposium
  on Code Generation and Optimization}.\hskip 1em plus 0.5em minus 0.4em\relax
  IEEE, 2016, pp. 182--193.

\bibitem{lam2019floatsmith}
M.~O. Lam, T.~Vanderbruggen, H.~Menon, and M.~Schordan, ``Tool integration for
  source-level mixed precision,'' in \emph{2019 IEEE/ACM 3rd International
  Workshop on Software Correctness for HPC Applications (Correctness)}, 2019,
  pp. 27--35.

\bibitem{clad}
V.~Vassilev, M.~Vassilev, A.~Penev, L.~Moneta, and V.~Ilieva, ``{Clad —
  Automatic Differentiation Using Clang and LLVM},'' \emph{Journal of Physics:
  Conference Series}, vol. 608, no.~1, p. 012055, 2015,
  [\href{http://stacks.iop.org/1742-6596/608/i=1/a=012055}{Link}].

\bibitem{8766229}
``Ieee standard for floating-point arithmetic,'' \emph{IEEE Std 754-2019
  (Revision of IEEE 754-2008)}, pp. 1--84, 2019.

\bibitem{betancourt2018geometric}
M.~Betancourt, ``{A geometric theory of higher-order automatic
  differentiation},'' \emph{arXiv preprint arXiv:1812.11592}, 2018.

\bibitem{griewank2008evaluating}
A.~Griewank and A.~Walther, \emph{{Evaluating derivatives: principles and
  techniques of algorithmic differentiation}}.\hskip 1em plus 0.5em minus
  0.4em\relax SIAM, 2008.

\bibitem{griewank1996algorithm}
A.~Griewank, D.~Juedes, and J.~Utke, ``{Algorithm 755: ADOL-C: A package for
  the automatic differentiation of algorithms written in C/C++},'' \emph{ACM
  Transactions on Mathematical Software (TOMS)}, vol.~22, no.~2, pp. 131--167,
  1996.

\bibitem{hogan2014fast}
R.~J. Hogan, ``{Fast reverse-mode automatic differentiation using expression
  templates in C++},'' \emph{ACM Transactions on Mathematical Software (TOMS)},
  vol.~40, no.~4, pp. 1--16, 2014.

\bibitem{sagebaum2018expression}
M.~Sagebaum, T.~Albring, and N.~R. Gauger, ``{Expression templates for primal
  value taping in the reverse mode of algorithmic differentiation},''
  \emph{Opt. Methods and Software}, vol.~33, no. 4-6, pp. 1207--1231, 2018.

\bibitem{hascoet2013tapenade}
L.~Hascoet and V.~Pascual, ``{The Tapenade automatic differentiation tool:
  principles, model, and specification},'' \emph{ACM Transactions on
  Mathematical Software (TOMS)}, vol.~39, no.~3, pp. 1--43, 2013.

\bibitem{enzymeNeurips}
W.~S. Moses and V.~Churavy, ``{Instead of Rewriting Foreign Code for ML,
  Automatically Synthesize Fast Gradients},'' in \emph{Advances in Neural
  Information Processing Systems 33}.\hskip 1em plus 0.5em minus 0.4em\relax
  Curran Associates, Inc., 2020.

\bibitem{princetontiger}
``{Princeton Research Computing},''
  [\href{https://researchcomputing.princeton.edu/systems/tiger}{Link}],
  accessed 2023-02-20.

\bibitem{llnlquartz}
``{Livermore Computing: HPC at LLNL},''
  [\href{https://hpc.llnl.gov/hardware/compute-platforms/quartz}{Link}],
  accessed 2023-02-20.

\bibitem{5306797}
S.~Che, M.~Boyer, J.~Meng, D.~Tarjan, J.~W. Sheaffer, S.-H. Lee, and
  K.~Skadron, ``Rodinia: A benchmark suite for heterogeneous computing,'' in
  \emph{2009 IEEE International Symposium on Workload Characterization
  (IISWC)}, 2009, pp. 44--54.

\bibitem{bienia11benchmarking}
C.~Bienia, ``Benchmarking modern multiprocessors,'' Ph.D. dissertation,
  Princeton University, January 2011.

\bibitem{fastapprox}
P.~Mineiro, ``Fastapprox,''
  [\href{https://code.google.com/archive/p/fastapprox/}{Google Archive Link}],
  2011.

\bibitem{pinnow2019typeforge}
N.~T. Pinnow, M.~Schordan, and T.~L. Vanderbrugger, ``Typeforge,'' Lawrence
  Livermore National Lab.(LLNL), Livermore, CA (United States), Tech. Rep.,
  2019.

\bibitem{brown2007profiling}
A.~W. Brown, P.~H. Kelly, and W.~Luk, ``Profiling floating point value ranges
  for reconfigurable implementation,'' in \emph{Proceedings of the 1st HiPEAC
  Workshop on Reconfigurable Computing}, 2007, pp. 6--16.

\bibitem{benz2012dynamic}
F.~Benz, A.~Hildebrandt, and S.~Hack, ``A dynamic program analysis to find
  floating-point accuracy problems,'' \emph{ACM SIGPLAN Notices}, vol.~47,
  no.~6, pp. 453--462, 2012.

\bibitem{lam2013dynamic}
M.~O. Lam, J.~K. Hollingsworth, and G.~Stewart, ``Dynamic floating-point
  cancellation detection,'' \emph{Parallel Computing}, vol.~39, no.~3, pp.
  146--155, 2013.

\bibitem{an2008fpinst}
D.~An, R.~Blue, M.~Lam, S.~Piper, and G.~Stoker, ``Fpinst: Floating point error
  analysis using dyninst,'' 2008.

\bibitem{das2020scalable}
A.~Das, I.~Briggs, G.~Gopalakrishnan, S.~Krishnamoorthy, and P.~Panchekha,
  ``Scalable yet rigorous floating-point error analysis,'' in \emph{SC20:
  International Conference for High Performance Computing, Networking, Storage
  and Analysis}.\hskip 1em plus 0.5em minus 0.4em\relax IEEE, 2020, pp. 1--14.

\bibitem{feliu2018abstract}
M.~A. Feli{\'u}, M.~Moscato, C.~A. Mu{\~n}oz \emph{et~al.}, ``An abstract
  interpretation framework for the round-off error analysis of floating-point
  programs,'' in \emph{International Conference on Verification, Model
  Checking, and Abstract Interpretation}.\hskip 1em plus 0.5em minus
  0.4em\relax Springer, 2018, pp. 516--537.

\bibitem{daumas2010certification}
M.~Daumas and G.~Melquiond, ``Certification of bounds on expressions involving
  rounded operators,'' \emph{ACM Transactions on Mathematical Software (TOMS)},
  vol.~37, no.~1, pp. 1--20, 2010.

\bibitem{lam2013automatically}
M.~O. Lam, J.~K. Hollingsworth, B.~R. de~Supinski, and M.~P. LeGendre,
  ``Automatically adapting programs for mixed-precision floating-point
  computation,'' in \emph{Proceedings of the 27th international ACM conference
  on International conference on supercomputing}, 2013, pp. 369--378.

\bibitem{guo2018exploiting}
H.~Guo and C.~Rubio-Gonz{\'a}lez, ``Exploiting community structure for
  floating-point precision tuning,'' in \emph{Proc. of the 27th ACM SIGSOFT
  Int. Symposium on Software Testing and Analysis}, 2018, pp. 333--343.

\bibitem{laguna2019gpumixer}
I.~Laguna, P.~C. Wood \emph{et~al.}, ``Gpumixer: Performance-driven
  floating-point tuning for gpu scientific applications,'' in \emph{Int. Conf.
  on High Performance Computing}.\hskip 1em plus 0.5em minus 0.4em\relax
  Springer, 2019, pp. 227--246.

\bibitem{iri1991history}
M.~Iri, ``History of automatic differentiation and rounding error estimation,''
  \emph{Andreas Griewank and George Corliss, editors}, pp. 3--16, 1991.

\bibitem{langlois2000revised}
P.~Langlois, ``A revised presentation of the cena method,'' Ph.D. dissertation,
  INRIA, 2000.

\bibitem{SaAlGauTOMS2019}
N.~G. M.~Sagebaum, T.~Albring, ``High-performance derivative computations using
  codipack,'' \emph{ACM Transactions on Mathematical Software (TOMS)}, vol.~45,
  no.~4, 2019.

\bibitem{quinlan2011rose}
D.~Quinlan and C.~Liao, ``The {ROSE} source-to-source compiler
  infrastructure,'' in \emph{Cetus users and compiler infrastructure workshop,
  in conjunction with PACT}, vol. 2011.\hskip 1em plus 0.5em minus 0.4em\relax
  Citeseer, 2011, p.~1.

\end{thebibliography}

\end{document}